\documentstyle[amscd]{amsart}
\newsymbol\varkappa 207B
\newsymbol\vartriangle 134D
\newsymbol\varnothing 203F
\newsymbol\restr 1316
\newsymbol\between 1347
\newsymbol\daleth 206B
\newsymbol\beth 2069
\newsymbol\twoheadrightarrow 1310
\newsymbol\Rsh 131F
\makeatletter
\renewcommand{\subsection}{\@startsection{subsection}{2}{\z@}%
{\baselineskip}{0.5\baselineskip}{\defaultfont\bf}}
\makeatother
\newtheorem{thm}{Theorem}[section]
\newtheorem{pro}[thm]{Proposition}
\newtheorem{lem}[thm]{Lemma}
\theoremstyle{definition}
\newtheorem{defn}{Definition}[section]
\numberwithin{equation}{section}
\def\dj{d\kern-.30em\raise1.25ex\vbox{\hrule width .3em height .03em}}
\def\Dj{D\kern-.75em\raise0.75ex\vbox{\hrule width .3em height .03em}
\kern.03em}
\newcommand{\rig}{\wp}
\newcommand{\hor}{\mbox{\family{euf}\shape{n}\selectfont hor}}
\newcommand{\ver}{\mbox{\family{euf}\shape{n}\selectfont ver}}
\newcommand{\con}{\mathrm{con}}
\newcommand{\adj}{\varpi}
\def\sstar{{\raise0.2ex\hbox{$\scriptstyle\star$}}}
\def\btwn{\gamma}
\def\amaph#1{*_{#1}}
\def\HZ{H\!Z}
\def\W{\frak{w}}
\def\map#1{#1_\star}
\def\1{\varnothing}
\def\3{\vartriangle}

\def\S{\Sigma}
\def\J{J_{\!\sigma}}
\def\IS{I(\Sigma)}
\def\IG{I^\sstar}
\def\AG{\cal{A}^{\!\mathrm{ad}}}
\def\r{\frak{r}}
\def\spl{\cal{X}}
\def\bla#1{$(${\it #1\/{}}$)$}
\def\Mor{\mathrm{Mor}}
\def\k{\kappa}
\def\im{\mathrm{im}}
\def\RepG{R(G)}
\def\Ch{\frak{ch}}
\def\hbim#1{\cal{F}_{#1}}
\def\bim#1{\cal{E}_{#1}}
\def\lhom{\mathrm{lhom}}
\def\rhom{\mathrm{rhom}}
\def\ghom{\mathrm{ghom}}
\def\hom{\mathrm{hom}}
\def\der{\mathrm{der}}
\def\rder{\mathrm{rder}}
\def\lder{\mathrm{lder}}
\def\WG{\Omega}
\def\WGi{\daleth}
\def\fWG{\widetilde{\adj}}
\def\hh{\frak{h}}
\def\rWG{\adj_\wedge}
\def\WGJ{\cal{K}}
\def\WGr{\Omega_*}
\def\WGri{\beth}
\def\tr{\mathrm{tr}}
\newcommand{\grten}{\mathbin{\widehat{\otimes}}}
\newcommand{\vh}{\mbox{\family{euf}\shape{n}\selectfont vh}}
\newcommand{\ad}{\mathrm{ad}}
\newcommand{\id}{\mathrm{id}}
\newcommand{\lie}{\mathrm{lie}}
\newcommand{\inv}{{i\!\hspace{0.8pt}n\!\hspace{0.6pt}v}}
\newcommand{\e}{\epsilon}
\newcommand{\Sum}{\displaystyle{\sum}}
\newcommand{\infM}{\frak{h}_M}
\newcommand{\hinfM}{H(\infM)}
\def\Gr{\frak{gr}}
\newcommand{\kI}{\Rsh}
\newcommand{\kW}{\amalg}
\begin{document}
\title[Quantum Principal Bundles]{Characteristic Classes\\
Of Quantum Principal Bundles}
\author{Mi\'co {\Dj}ur{\Dj}evi\'c}
\address{Instituto de Matematicas, UNAM, Area de la Investigacion
Cientifica, Circuito Exterior, Ciudad Universitaria, M\'exico DF, CP
04510, MEXICO}
\maketitle
\begin{abstract}
A noncommutative-geometric generalization of classical Weil theory of
characteristic classes is presented, in the conceptual framework of
quantum principal bundles. A particular care is given to the case when
the bundle does not admit regular connections. A cohomological
description of the domain of the Weil homomorphism is given. Relations
between universal characteristic classes for the regular and
the general case are analyzed. In analogy with classical
geometry, a natural spectral sequence is introduced and investigated.
The appropriate counterpart of the Chern character is constructed, for
structures admitting regular connections. Illustrative examples and
constructions are presented.
\end{abstract}
\renewcommand{\thepage}{}
\tableofcontents
\section{Introduction}

The theory of characteristic classes plays a fundamental role in
understanding the internal structure of classical spaces. In classical
differential geometry, the Weil approach \cite{KN} to characteristic
classes is the most natural. Characterisctic classes are labeled by invariant
polynomials on the Lie algebra of the structure group. A connection with
the purely topological approach is given by establishing the
isomorphism between the algebra of invariant polynomials and the
cohomology algebra of the classifying space of the structure Lie group.

The aim of this study is to incorporate classical Weil theory of
characteristic classes into a general quantum framework, in accordance
with basic principles of non-commutative differential geometry
\cite{C1,C2}.
Considerations of this study are logically based on a general theory
of quantum principal bundles, developed in \cite{D1,D2}.

\renewcommand{\thepage}{\arabic{page}}
Let us outline the contents of the paper. In the next section the main
facts from the theory of quantum principal
bundles are collected. This includes a brief exposition of differential
calculus, and the formalism of connections, together with the
description of main associated algebras. In particular, general
constructions of operators of covariant derivative, horizontal
projection and curvature are sketched. Finally, the structure
of the algebra of horizontal forms is analyzed, from the point of view
of the representation theory of the structure group.

General theory of characteristic classes is presented in Section 3. As
first, a construction of the Weil homomorphism given in \cite{D2} is
sketched. This works fine when the bundle admits regular and multiplicative
connections. However it turns out that essentially the same
construction can be performed in a more general context, when
the algebra of horizontal forms and the appropriate space of
``vertical'' forms are mutually transversal. Next, we shall discuss a
question of renormalizing the calculus on the bundle so that an
arbitrarily given connection becomes regular, in terms of the
renormalized calculus. Such a construction always works, however the
resulting calculus may be too simple, if the initial connection is
sufficiently irregular.

We shall then introduce and analyze a natural spectral sequence, in
analogy with classical geometry. This spectral sequence is associated to
the filtration of the calculus on the bundle, induced by the pull back
of the right action map. The corresponding first two cohomology
algebras will be analyzed in details.

The last and the main theme of Section 3 is a presentation of a universal
approach to characteristic classes, which gives a cohomological
interpretation of the domain of the Weil homomorphism,
incorporating certain elements of classical classifying space
interpretation of universal characteristic classes (as cohomology
classes of the corresponding classifying space).

The construction
works for arbitrary quantum principal bundles and connections. It also
admits a natural ``projection'' to the level of bundles admitting
regular connections. Interestingly, although the two levels of the
theory (regular and general) are formulated in essentially the same way,
the corresponding characteristic classes radically differ in
sufficiently ``irregular'' situations.

For example,  at the level of general bundles there exist generally
nontrivial classes in {\it odd} dimensions, in contrast to the regular
case where all classes are expressible in terms of the curvature map.
Another important difference between two levels is in the
existence of the Chern-Simons forms. At the general level, all
characteristic classes vanish, as cohomology classes of the bundle. In
contast to this, there exist regular characteristic classes, nontrivial
as classes on the bundle. In particular, such classes have no
analogs at the level of general differential structures. However, it
turns out that under certain regularity assumptions on the calculus over
the structure group, universal classes for the regular and the general
case coincide. This essentially simplifies the work with bundles without
regular connections.

The basic theme of Section 4 is the construction of the quantum Chern
character. It is essential for the construction that the bundle admits
regular connections. The first step is to construct a pairing between
vector bundles associated to a given quantum principal bundle $P$, and
differential forms from the graded centre $Z(M)$ of the algebra
$\Omega(M)$ representing the calculus on the base manifold $M$.
Associated vector bundles will be represented by intertwiner
$\Omega(M)$-bimodules, between representations of the structure
group $G$ and the right action of $G$ on horizontal forms
(so that the elements of these bimodules are interpretable as
vector-bundle valued forms on $M$). The above mentioned pairing
will be constructed with the help of a canonical
``trace'' on the corresponding vector bundle endomorphism algebras, and
the covariant derivative of an arbitrary regular connection.

The next step in the construction is to pass to the corresponding
cohomology algebra $\HZ(M)$. This rules out the connection-dependence of
the pairing. Finally, starting from intertwining bimodules (with natural
operations of addition and tensor product over $\Omega(M)$), we
introduce the corresponding $K$-ring $K(M,P)$. In such a way we obtain
the quantum Chern character, as a right $K(M,P)$-module structure
$\Ch_M\colon \HZ(M)\times K(M,P)\rightarrow \HZ(M)$. Finally, we find an
explicit expression for $\Ch_M$, in terms of the curvature map.

At the end of Section 4 we present a construction of quantum Chern
classes, generalizing the classical ``determinantal''
definition. However, an interesting quantum phenomena appears--generally
a given vector bundle will generate classes in arbitrary high
dimensions.

In Section 5 properties of quantum $U(1)$-bundles are studied. A
complete structural analysis of these bundles is performed, including
the classification of differential calculi on them. Here the general
formalism essentially simplifies, as far as we assume that the calculus
on the structure group is $1$-dimensional. In particular, this rules
out the question about ``embedded differential'' maps, and  uniquely fixes the
``horizontal-vertical'' decomposition. The algebra of characteristic
classes is generated by a single element, the Euler class, given by the
curvature form.

We shall also introduce the Euler class for certain vector bundles
associated to general quantum principal bundles admitting regular
connections and appropriate braided volume forms, following the
classical analogy.

Section 6 is devoted to the study of characteristic classes of locally
trivial bundles  over classical smooth manifolds \cite{D1}. It turns
out that characteristic classes of such a quantum principal bundle $P$ can be
viewed as standard characteristic classes of its classical part $P_{cl}$
(the set of classical points of $P$). We shall also describe an
interesting class of line bundles which are not locally trivializable,
and discuss their characteristic classes.

In Section 7 we present a construction of a natural differential
calculus over quantum frame bundles. These structures generally give
examples with non-trivial characteristic classes. However all
computations can be performed completely because of the simple internal
structure of the calculus.

Finally, in the last section some concluding observations and
additional examples are collected. We
end this introduction with several remarks concerning the relevant
quantum group entities. Conceptually, we follow \cite{W1,W2} in a slightly
different notation.

Let us consider a compact matrix quantum group $G$, represented
by a Hopf *-algebra $\cal{A}$, interpreted as consisting
of polynomial functions on $G$. Let
$\phi,\e$ and $\k$ be the coproduct, counit and the antipode map
respectively. Let $\ad\colon\cal{A}\rightarrow\cal{A}\otimes\cal{A}$ be
the adjoint action of $G$ on itself.

     Let $\Gamma$ be a bicovariant *-calculus \cite{W2}
over $G$ and $\Gamma^{\wedge}$
its universal differential envelope (\cite{D1}--Appendix B).
Explicitly,
$$\Gamma^\wedge=\Gamma^\otimes/S^\wedge$$
where $\Gamma^\otimes$ is the tensor bundle algebra,
and $S^\wedge\subseteq\Gamma^\otimes$ is the ideal generated by elements
of the form
$$Q=\sum_k da_kdb_k\qquad\sum_k a_kdb_k=0.$$

The coproduct map $\phi\colon \cal{A}\rightarrow\cal{A}\otimes
\cal{A}$ is uniquely extendible to the homomorphism
$\widehat{\phi}\colon \Gamma^{\wedge}\rightarrow\Gamma^{\wedge}
\grten\Gamma^{\wedge}$ of graded-differential algebras. The restriction
of this map on $\Gamma$ is given by
\begin{equation*}
\widehat{\phi}(\vartheta)=\ell_\Gamma(\vartheta)+\rig_\Gamma(\vartheta),
\end{equation*}
where $\ell_\Gamma$ and $\rig_\Gamma$ are the corresponding left
and right action of $G$ on $\Gamma$.

Let $\Gamma_{\inv}$  be the left-invariant  part of $\Gamma$. We shall
denote by $\cal{R}\subseteq \ker(\e)$ the  right
$\cal{A}$-ideal which determines $\Gamma$. There exists a
natural surjection $\pi\colon \cal{A}\rightarrow\Gamma_{\inv}$ given by
\begin{equation*}
\pi(a)=\k(a^{(1)})da^{(2)}.
\end{equation*}
We have $\ker(\pi)=\Bbb{C}\oplus{\cal R}$, and hence there exists a
natural isomorphism
$$\Gamma_{\inv}\leftrightarrow\ker(\e)/\cal{R}.$$

The product in $\cal{A}$ induces a right $\cal{A}$-module structure on
$\Gamma_{\inv}$. We shall denote by $\circ$ this structure. Explicitly,
\begin{equation*}
\pi (a)\circ  b=\pi \bigl(ab-\e(a)b\bigr),
\end{equation*}
for each $a,b\in\cal{A}$.

Throughout this paper we shall deal with various bicovariant *-algebras
$\Gamma^\star$ over $\Gamma$. The symbol $\Gamma^\star_{\inv}$ will then
refer to the corresponding left-invariant subalgebra. Here ${}^\star$
symbolically characterizes the algebra---some other
entities naturally related to $\Gamma^\star$ will be also endowed
with the superscript ${}^\star$.

We have
$$\Gamma_{\inv}^\wedge=\Gamma_{\inv}^\otimes/S^\wedge_{\inv}$$
where $S^\wedge_{\inv}$ is the left-invariant part of $S^\wedge$. As an
ideal in $\Gamma_{\inv}^\otimes$, it is generated by elements of the
form
$$ q=\pi(a^{(1)})\otimes\pi(a^{(2)})$$
where $a\in\cal{R}$.

The $\circ$ can be naturally extended from $\Gamma_{\inv}$ to
$\Gamma_{\inv}^\star$, with the help of the formula
$$ \vartheta\circ  a=\k(a^{(1)})\vartheta a^{(2)}, $$
so that we have
\begin{align*}
1\circ a&=\e(a)1\\
(\vartheta\eta)\circ a&=(\vartheta\circ a^{(1)})(\eta\circ
a^{(2)}).
\end{align*}

The adjoint action $\adj \colon \Gamma_{\inv}
\rightarrow\Gamma_{\inv}\otimes \cal{A}$
of $G$ on $\Gamma_{\inv}$ is explicitly given by
$$ \adj \pi=(\pi\otimes  \id) \ad,$$
it is naturally extendible to algebras $\Gamma_{\inv}^\sstar$.

Let $\widehat{\adj}\colon\Gamma_{\inv}^\wedge\rightarrow
\Gamma_{\inv}^\wedge\grten\Gamma^\wedge$ be the graded-differential
extension of $\adj$. We have
$$
\widehat{\adj} (\vartheta)=1\otimes \vartheta
+\adj (\vartheta),
$$
for each $\vartheta\in\Gamma_{\inv}$.

Finally, the following identities hold
\begin{gather*}
d(\vartheta\circ  a)=d(\vartheta)\circ  a-\pi(a^{(1)})(\vartheta\circ
a^{(2)})+(-1)^{\partial\vartheta}(\vartheta\circ
a^{(1)})\pi(a^{(2)})\\
d\pi(a)=-\pi(a^{(1)})\pi(a^{(2)})\\
(\vartheta\circ  a)^*=\vartheta^*\circ  \k(a)^*\qquad \pi(a)^*
=-\pi\left[\k(a)^*\right].
\end{gather*}
All the above listed properties will be commonly used in various
technical considerations.

\section{General Conceptual Framework}

We pass to a brief exposition of general theory of quantum principal
bundles, paying particular attention to differential calculus and the
formalism of connections.

Let us consider a quantum space $M$, represented by a *-algebra
$\cal{V}$. Let $P=(\cal{B},i,F)$ be a quantum principal $G$-bundle over
$M$. Here $\cal{B}$  is a *-algebra representing $P$ as a quantum space,
while $i\colon\cal{V}\rightarrow\cal{B}$ and
$F\colon\cal{B}\rightarrow\cal{B}\otimes\cal{A}$ are unital
*-homomorphisms playing the role of the dualized projection of $P$
on $M$ and the right action of $G$ on $P$ respectively.

Let us assume that a differential calculus on the structure quantum
group $G$ is based on the universal envelope $\Gamma^\wedge$ of a given
first-order bicovariant *-calculus $\Gamma$.
Let $\Omega(P)$ be a graded-differential *-algebra
representing the calculus on $P$. As first, this means that
$\Omega^0(P)=\cal{B}$, and
that $\cal{B}$ generates the differential algebra $\Omega(P)$.
Secondly the right action $F$ is extendible (necessarily uniquely) to
a graded-differential *-homomorphism $\widehat{F}\colon\Omega(P)
\rightarrow\Omega(P)\grten\Gamma^\wedge$.

The right action
$F^\wedge\colon\Omega(P)\rightarrow\Omega(P)\otimes\cal{A}$ of $G$ on
$\Omega(P)$ is given by
$$F^\wedge=(\id\otimes p_*)\widehat{F},$$
where $p_*\colon\Gamma^\wedge\rightarrow\cal{A}$ is the projection map.

Let us consider a connection $\omega\colon\Gamma_{\inv}
\rightarrow\Omega(P)$. This means that $\omega$ is a hermitian intertwiner
between $\adj\colon\Gamma_{\inv}\rightarrow\Gamma_{\inv}
\otimes\cal{A}$ and the right action $F^\wedge\colon\Omega(P)\rightarrow
\Omega(P)\otimes\cal{A}$, satisfying
$$\pi_v\omega(\vartheta)=1\otimes\vartheta$$
where $\pi_v\colon\Omega(P)\rightarrow\ver(P)$ is the
verticalization homomorphism, and $\ver(P)$ is the
graded-differential *-algebra representing ``verticalized'' differential
forms on the bundle (at the level of graded vector spaces, $\ver(P)=\cal{B}
\otimes\Gamma_{\inv}^\wedge$). The bundle $P$ possesses at least
one connection. The set $\con(P)$ of all connections on $P$ is a
real affine space, in a natural manner.

Connections can be equivalently described as hermitian first-order
linear maps $\omega\colon\Gamma_{\inv}\rightarrow\Omega(P)$ satisfying
\begin{equation}\label{con}
\widehat{F}\omega(\vartheta)=\sum_k\omega(\vartheta_k)\otimes
c_k+1\otimes\vartheta,
\end{equation}
where $\adj(\vartheta)=\Sum_k\vartheta_k\otimes c_k$.

The graded *-subalgebra
of horizontal forms is defined by
$$\hor(P)=\widehat{F}^{-1}\Bigl\{\Omega(P)\otimes\cal{A}\Bigr\}.$$
It turns out that $\hor(P)$ is $F^\wedge$-invariant.

Let $\Omega(M)\subseteq\Omega(P)$ be the graded-differential *-subalgebra
representing differential forms on the base manifold $M$. The
elements of $\Omega(M)$ are characterized as $\widehat{F}$-invariant
differential forms on $P$. Equivalently, $\Omega(M)$ is a
$F^\wedge$-fixed point subalgebra of $\hor(P)$. In some formulas we shall
use a special symbol $d_M\colon\Omega(M)\rightarrow\Omega(M)$ for
the corresponding differential.

By definition, $\omega$ is {\it regular} if the following
commutation relation holds
\begin{equation}
\omega(\vartheta)\varphi=(-1)^{\partial\varphi}\sum_k\varphi_k\omega(\vartheta
\circ c_k)
\end{equation}
where $\varphi$ is an arbitrary horizontal form, and
$\Sum_k\varphi_k\otimes c_k=F^\wedge(\varphi)$.
In particular, regular connections graded-commute with the
elements of $\Omega(M)$.

A connection $\omega$ is called {\it multiplicative} iff it is extendible
(necessarily uniquely) to a (unital) *-homomorphism $\omega^\wedge\colon
\Gamma^\wedge_{\inv}\rightarrow\Omega(P)$.

It is important to observe that multiplicativity and regularity
properties essentially depend on the complete differential calculus on
the bundle. In general, $\Omega(P)$ may be such that there will be no regular
nor multiplicative connections. Furthermore, regular connections, if
exist, are not necessarily multiplicative. However, if one regular
connection is multiplicative, then all regular connections will be
multiplicative, too. If regular connections are not multiplicative then
it is possible to ``renormalize'' the calculus on the bundle, by
factorizing $\Omega(P)$ through an appropriate graded-differential
*-ideal, measuring a lack of multiplicativity of regular
connections. Such a factorization does not change the first-order
calculus. In terms of the renormalized calculus, regular connections are
multiplicative.

The map $\omega^\wedge\colon\Gamma_{\inv}^\wedge\rightarrow\Omega(P)$
can be introduced for an arbitrary connection, by the formula
$$\omega^\wedge=\omega^\otimes\iota$$
where $\iota\colon\Gamma_{\inv}^\wedge\rightarrow\Gamma_{\inv}^\otimes$ is an
arbitrary grade-preserving hermitian section of the
factorization map intertwining the adjoint actions of $G$, and
$\omega^\otimes\colon\Gamma_{\inv}^\otimes\rightarrow\Omega(P)$ is
the unital multiplicative extension of $\omega$.

The quantum counterparts of horizontal projection, covariant derivative
and curvature can be constructed as follows. Let $\vh(P)$ be the graded
*-algebra of ``vertically-horizontally'' decomposed forms. Explicitly
$\vh(P)=\hor(P)\otimes\Gamma_{\inv}^\wedge$, at the level of graded
vector spaces. The *-algebra structure is specified by
\begin{align*}
(\psi\otimes\eta)(\varphi\otimes\vartheta)&=
(-1)^{\partial\eta\partial\varphi}\sum_k\psi\varphi_k\otimes
(\eta\circ c_k)\vartheta\\
(\varphi\otimes\vartheta)^*&=\sum_k\varphi_k^*\otimes(\vartheta^*\circ
c_k^*).
\end{align*}
Let us consider the ``decomposition map''
$m_\omega\colon\vh(P)\rightarrow\Omega(P)$, given by
$$ m_\omega(\varphi\otimes\vartheta)=\varphi\omega^\wedge(\vartheta).$$
This map is bijective, and intertwines the corresponding
right actions of $G$.
Moreover, if $\omega$ is regular and multiplicative then $m_\omega$ is a
*-algebra isomorphism.

The horizontal projection operator $h_\omega\colon\Omega(P)\rightarrow
\hor(P)$ is given by
$$h_\omega=(\id\otimes p_*)m_\omega^{-1}.$$
In a direct analogy with classical differential geometry, the covariant
derivative operator can be defined by
$$D_\omega=h_\omega d.$$
The operators $D_\omega$ and $h_\omega$ are right-covariant, in a
natural manner. Finally, the curvature
$R_\omega\colon\Gamma_{\inv}\rightarrow \Omega(P)$
can be defined as the composition
$$R_\omega=D_\omega\omega.$$
This is a hermitian tensorial 2-form. The following identity holds
$$R_\omega=d\omega-\langle\omega,\omega\rangle.$$
Here $\langle\,\rangle$ are the
brackets associated to the embedded differential
$\delta\colon\Gamma_{\inv}\rightarrow\Gamma_{\inv}^{\otimes 2}$,
where $\delta(\vartheta)=\iota d(\vartheta)$. The above identity
generalizes the classical Structure equation.

It is worth noticing that $R_\omega$ depends, besides on $\omega$, also
on the map $\delta$. In a large part of the formalism, the section
$\iota$ figures only via the corresponding $\delta$. The map $\delta$ is
always related with the transposed commutator \cite{W2}
map $c^\top=(\id\otimes\pi)\adj$ in the following way
\begin{equation}\label{sd=d+c}
\sigma\delta=\delta+c^\top,
\end{equation}
where $\sigma\colon\Gamma_{\inv}^{\otimes 2}\rightarrow\Gamma_{\inv}^{
\otimes 2}$ is the (left-invariant restriction of the)
canonical braid operator \cite{W2}. It is explicitly given
by
\begin{equation}
\sigma(\eta\otimes\vartheta)=\sum_k\vartheta_k\otimes(\eta\circ c_k),
\end{equation}
where $\adj(\vartheta)=\Sum_k\vartheta_k\otimes c_k$.

The map $c^\top$ is extendible to
$c^\top\colon\Gamma_{\inv}^\sstar\rightarrow\Gamma_{\inv}^\sstar
\otimes\Gamma_{\inv}$, by the same formula.
The elements of $S_{\inv}^{\wedge2}$ are $\sigma$-invariant, and hence
$\sigma$ projects to $\sigma\colon\Gamma_{\inv}^{\wedge
2}\rightarrow\Gamma_{\inv}^{\wedge2}$.

Let $\IG\subseteq\Gamma_{\inv}^\otimes$ be the graded *-subalgebra
consisting of elements invariant under the adjoint action $\adj\colon
\Gamma^\otimes_{\inv}\rightarrow\Gamma_{\inv}^\otimes\otimes\cal{A}$. The
formula
\begin{equation}\label{weil*}
W^\omega(\vartheta)=R_\omega^\otimes(\vartheta)
\end{equation}
defines a *-homomorphism $W^\omega\colon\IG\rightarrow\Omega(M)$.
Here $R_\omega^\otimes\colon\Gamma_{\inv}^\otimes\rightarrow\Omega(P)$ is
the corresponding unital multiplicative extension.

We shall denote by $Z(M)$ the graded centre of
$\Omega(M)$. It is a graded-differential *-subalgebra of $\Omega(M)$.

If $\omega$ is regular then the following identity holds
$$ R_\omega(\vartheta)\varphi=\sum_k\varphi_kR_\omega(\vartheta\circ
c_k),$$
for each $\varphi\in\hor(P)$.
In particular, the curvature $R_\omega$ commutes with all elements of
$\Omega(M)$, and hence the image of $W^\omega$ is contained in $Z(M)$.

For the end of this section, we consider the internal
structure of the algebras $\cal{B}$ and $\hor(P)$, from the point
of view of the representation theory \cite{W2} of
compact matrix quantum groups. The symbol $\otimes_M$ will be used for
both tensor products over $\cal{V}$ and $\Omega(M)$, depending on the
context.

Let $\RepG$ be a complete monoidal category of finite-dimensional
unitary representations of $G$. Each $u\in\RepG$ is understandable
as a map $u\colon H_u\rightarrow H_u\otimes\cal{A}$, where $H_u$ is the
corresponding carrier unitary space (we shall assume that $H_u$ belong
to a fixed model set of unitary spaces). Morphisms of the
category $\RepG$ are intertwiners between
representations. We shall use the symbols
$\oplus,\times$ for the sum and the product in $\RepG$. We have
$$ H_{u\oplus v}=H_u\oplus H_v\qquad H_{u\times v}=H_u\otimes H_v.$$
The symbol $\1$ will refer to the {\it trivial}
representation of $G$ acting in the space $H_{\!\1}=\Bbb{C}$.

For each $u\in\RepG$
let $\bim{u}=\Mor(u,F)$ be the space of intertwiners between $u$
and $F\colon\cal{B}\rightarrow\cal{B}\otimes\cal{A}$. The spaces
$\bim{u}$ are $\cal{V}$-bimodules, in a natural manner.
Furthermore $\bim{\1}=\cal{V}$, in a natural manner.

For each $u,v\in\RepG$, the following natural bimodule isomorphism
holds:
\begin{equation}
\bim{u\times v}\leftrightarrow\bim{u}\otimes_M\!\bim{v}.
\end{equation}
The above
isomorphism is induced by the product in $\cal{B}$. Explicitly
$$ \varphi\otimes\psi\colon x\otimes y\mapsto\varphi(x)\psi(y). $$

Let us denote by $u^{c}=\bar{u}$ the
contragradient representation, identified with the
conjugate representation, acting in $H_{u^{\!c}}=H_u^*$.
Let $C_u\colon H_u\rightarrow H_u$ be the canonical intertwiner between
$u$ and the second contragradient $u^{cc}$.
This map is strictly positive, and normalized such
that $\tr(C_u)=\tr(C_u^{-1})$. Furthermore,
$$ C_u\otimes C_v=C_{u\times v}\qquad \varphi C_u=C_v\varphi$$
for each $u,v\in\RepG$ and $\varphi\in\Mor(u,v)$. The appropriate scalar
product in $H_u^*$ is given by
$$ (f,g)=(j_u^{-1}(g),C_uj_u^{-1}(f))$$
where $j_u\colon H_u\rightarrow H_u^*$ is the antilinear map induced by
the scalar product in $H_u$.

There exists a natural bimodule anti-isomorphism $\amaph{u}\colon
\bim{u}\rightarrow\bim{\bar{u}}$, determined by the diagram
\begin{equation}
\begin{CD}
H_{u} @>{\mbox{$\phantom{\amaph{u}}\varphi$}}>> \cal{B}\\
@V{\mbox{$j_u$}}VV @VV{\mbox{$*$}}V\\
H_{u}^* @>>{\mbox{$\amaph{u}\varphi$}}> \cal{B}
\end{CD}
\end{equation}
In other words, we can naturally identify $\bim{\bar{u}}$ and the
corresponding conjugate bimodule $\bim{u}^*$.

Let $\rhom(\bim{u},\bim{v})$ be the space of right $\cal{V}$-linear
maps $A\colon\bim{u}\rightarrow\bim{v}$. This is a
$\cal{V}$-bimodule, in a natural manner. Similarly, we shall denote by
$\lhom(\bim{u},\bim{v})$ the space of corresponding
left $\cal{V}$-module homomorphisms. Finally, let us denote by
$\hom=\lhom\cap\rhom$ the spaces of bimodule homomorphisms.
For $u=v$ all introduced spaces are algebras.

Each map $f\in\Mor(u,v)$ induces a bimodule
homomorphism $\map{f}\colon\bim{v}\rightarrow\bim{u}$, given by the
composition of intertwiners.

{}From the point of view of our
considerations, the most interesting intertwiners are given by the
standard contraction map $\btwn^u\colon H_u^*\otimes H_u\rightarrow\Bbb{C}$,
and the natural inclusion $I_u\colon\Bbb{C}\rightarrow H_u\otimes
H_u^*=L(H_u)$. These intertwiners induce bimodule homomorphisms
$\map{\btwn^u}\colon\cal{V}\rightarrow\bim{u}^*\otimes_M\!\bim{u}$, and
$\langle\rangle_u^+\colon\bim{u}\otimes_M\!\bim{u}^*\rightarrow\cal{V}$
respectively. The map $\map{\btwn^u}$ is completely determined by
$$
\map{\btwn^u}(1)=\sum_k\nu_k\otimes\mu_k.
$$
By definition, the following identity holds
\begin{equation}
\sum_k\nu_k(e_i^*)\mu_k(e_j)=\delta_{ij}1,
\end{equation}
where $\{e_i\}$ is an arbitrary orthonormal basis in $H_u$ and
$\{e_i^*\}$ is the corresponding biorthogonal basis in $H_u^*$.

The map $\langle\rangle_u^+$ is naturally understandable as a pairing
between $\bim{u}$ and $\bim{u}^*$. Explicitly,
\begin{equation}
\langle\varphi,\psi\rangle_u^+=\sum_i\varphi(e_i)\psi(e_i^*).
\end{equation}

Similarly, it is possible to reverse the roles of $u$ and $\bar{u}$. In
particular, we can construct the canonical pairing $\langle\rangle_u^-
\colon\bim{u}^*\otimes_M\!\bim{u}\rightarrow\cal{V}$. It is given by
\begin{equation}
\langle\psi,\varphi\rangle_u^-=\sum_{ij}(C_u^{-1})_{ji}\psi(e_i^*)
\varphi(e_j).
\end{equation}
The operator $C_u^{-1}$ in the above expression comes from the
compensation of the appearance of the second conjugate $u^{cc}$
in place of $u$.

The structure of the bundle
$P$ is expressible \cite{D3} in terms of the system of
bimodules $\bim{u}$, together with maps of
the form $\amaph{u}$ and $\map{f}$. In particular,
the following natural decomposition holds
\begin{equation}
\cal{B}=\sideset{}{^\oplus}\sum_{\alpha\in\cal{T}}\cal{B}^\alpha,
\end{equation}
where $\cal{T}$ is the complete set
of mutually inequivalent irreducible unitary representations of $G$,
and $\cal{B}^\alpha$ are multiple irreducible
$\cal{V}$-submodules with respect to the decomposition of $F$.
These modules are further naturally decomposable as follows
\begin{equation}
\cal{B}^\alpha=\bim{\alpha}\otimes H_\alpha\qquad
\varphi(x)\leftrightarrow\varphi\otimes x.
\end{equation}

A similar consideration can be applied to the algebra $\hor(P)$. This
leads to the intertwiner $\Omega(M)$-bimodules
$\hbim{u}=\Mor(u,F^\wedge)$. We shall use the same symbols for the
characteristic spaces and operators related to these bimodules. The
spaces $\hbim{u}$ are naturally graded, and obviously
$\hbim{u}^0=\bim{u}$. Moreover, it turns out that right/left product
maps in $\hbim{u}$ induce canonical decompositions
\begin{equation}\label{WM-M}
\bim{u}\otimes_M\!\Omega(M)\leftrightarrow\hbim{u}\leftrightarrow
\Omega(M)\otimes_M\!\bim{u}.
\end{equation}
In particular $\lhom(\bim{u},\bim{v})$ is realized as a subalgebra of
$\lhom(\hbim{u},\hbim{v})$, consisting of grade-preserving
transformations.

The structure of spaces $\hbim{u}$ can be expressed in terms of spaces
$\bim{u}$. Composing the above two identifications we obtain canonical
flip-over maps $$\sigma_u\colon\bim{u}\otimes_M\!\Omega(M)
\rightarrow\Omega(M)\otimes_M\!\bim{u}.$$ By definition, these maps are
grade-preserving and act as the identity on $\bim{u}$. The
$\Omega(M)$-bimodule structure on $\hbim{u}$ is expressed by the
following compatibility condition
\begin{equation}
\sigma_u(\id\otimes m_M)=(m_M\otimes\id)
(\id\otimes\sigma_u)(\sigma_u\otimes\id),
\end{equation}
where $m_M$ is the product in $\Omega(M)$. Intertwiner homomorphisms
$\map{f}$ and conjugation maps $\amaph{u}$ satisfy the following
diagrams
\begin{equation*}
\begin{CD} \bim{v}\otimes_M\!\Omega(M) @>{\mbox{$\sigma_v$}}>>
\Omega(M)\otimes_M\!\bim{v}\\
@V{\mbox{$\map{f}\otimes\id$}}VV @VV{\mbox{$\id\otimes\map{f}$}}V\\
\bim{u}\otimes_M\!\Omega(M) @>>{\mbox{$\sigma_u$}}>
\Omega(M)\otimes_M\!\bim{u}
\end{CD}\qquad
\begin{CD} \bim{u}\otimes_M\!\Omega(M) @>{\mbox{$\sigma_u$}}>>
\Omega(M)\otimes_M\!\bim{u}\\
@V{\mbox{$\amaph{u}$}}VV  @VV{\mbox{$\amaph{u}$}}V\\
\Omega(M)\otimes_M\!\bim{\bar{u}} @<<{\mbox{$\sigma_{\bar{u}}$}}<
\bim{\bar{u}}\otimes_M\!\Omega(M)
\end{CD}
\end{equation*}
The first diagram actually  characterises elements of
$\hom(\bim{v},\bim{u})$ that are extendible to corresponding
$\Omega(M)$-bimodule
homomorphisms (the left and right $\Omega(M)$-linear extensions
coincide).

The algebra of horizontal forms can be decomposed in the following way
\begin{equation}
\begin{gathered}
\hor(P)=\sideset{}{^\oplus}\sum_{\alpha\in\cal{T}}\cal{H}^\alpha(P)\qquad
\cal{H}^\alpha(P)=\hbim{\alpha}\otimes H_\alpha\\
\hor(P)\leftrightarrow\Omega(M)\otimes_M\!\cal{B}
\leftrightarrow\cal{B}\otimes_M\!\Omega(M).
\end{gathered}
\end{equation}

Let us observe that the elements
of $\lhom(\hbim{u},\hbim{v})$ are also naturally interpretable as maps
$A$ satisfying
\begin{equation}\label{gr-lin}
A(\alpha\psi)=(-1)^{\partial\alpha\partial \!A}\alpha A(\psi).
\end{equation}

In what follows we shall commonly pass from one interpretation to
another. We shall denote by $\ghom(\hbim{u},\hbim{v})$ the
space of elements $A\in\rhom(\hbim{u},\hbim{v})$ satisfying also
\eqref{gr-lin}.

Let us now consider a regular connection $\omega$ on $P$.
The covariant derivative map $D_\omega$ is a right-covariant
hermitian first-order derivation on $\hor(P)$.
The composition of $D_\omega$ with elements
of $\hbim{u}$ induces first-order maps
$D_{\omega,u}\colon\hbim{u}\rightarrow\hbim{u}$.

The graded Leibniz rule for $D_\omega$ and the
fact that $(D_\omega-d_M)\bigl(\Omega(M)\bigr)=\{0\}$ imply
$$ D_{\omega,u\times v}(\varphi\otimes\psi)=D_{\omega,u}
(\varphi)\otimes\psi
+(-1)^{\partial\varphi}\varphi\otimes D_{\omega,v}(\psi)\qquad
D_{\omega,\1}=d_M.$$
Moreover, the following compatibility properties hold
$$ D_{\omega,\bar{u}}\amaph{u}=\amaph{u}D_{\omega,u}\qquad
D_{\omega,u}\map{f}=\map{f}D_{\omega,v}.$$

\section{Quantum Characteristic Classes}

In this section the abstract theory of characteristic classes will be
presented. We shall start from a particular case, when the bundle admits
regular connections (and the caclulus is ``renormalized'' such that
regular connections are multiplicative). Then the construction of the
Weil homomorphism can be performed in a direct analogy with the
classical case. Then we pass to a more general geometrical framework,
which however still allows to perform the constructions as in the
regular case. Finally, we present a general cohomological
approach to characteristic classes, in the framework of which ``universal''
characteristic classes are defined. We shall consider separately the
general and the regular case. These universal classes naturally
form the domain of the Weil homomorphism. A particular attention will be
payed to the study of interrelations between two levels of generality.
The formalism is applicable to arbitrary bundles.

\subsection{The Regular Case}
Let us assume that the bundle admits regular (and multiplicative)
connections. Then it is possible to prove, straightforwardly
generalizing the classical exposition \cite{KN}, that
$dW^{\omega}(\vartheta)=0$ for each $\vartheta \in \IG$
and $\omega\in\r(P)$. Moreover \cite{D2}, the cohomological
class of $W^{\omega}(\vartheta)$ in $\Omega (M)$ is
independent of the choice of $\omega\in\r(P)$.
\begin{lem}
Actually a finer result
holds---the above cohomological class is well defined in $Z(M)$.\qed
\end{lem}

In such a way we have constructed an intrinsic unital *-homomorphism
$$W\colon\IG\rightarrow \HZ(M)\qquad
W(\vartheta)=[W^{\omega}(\vartheta)],$$
where $\HZ(M)$ is the corresponding cohomology algebra. This map is a
quantum counterpart of the Weil homomorphism. We shall denote by
$\W(M)\subseteq \HZ(M)$ the image of $W$.

The quantum Weil homomorphism can be further factorized through
the ideal $\J$  generated  by the space
$\im(I-\sigma)\subseteq \Gamma_{\inv}^{\otimes 2}$. This follows
from the commutation identity
$$ R_\omega(\eta)R_\omega(\vartheta)=\sum_k
R_\omega(\vartheta_k)R_\omega(\eta\circ c_k) $$
where $\adj(\vartheta)=\Sum_k\vartheta_k\otimes c_k$.
The elements of the factoralgebra
$\S =\Gamma_{\inv}^\otimes/\J$
are interpretable as ``polynoms'' over the ``Lie  algebra''  of  $G$.
The  adjoint action $\adj^{\otimes}$  is naturally projectable on
$\S$.  Let  $\IS\subseteq \S$
be the  subalgebra  of  elements  invariant  under  the  projected
action (invariant polynomials). We have
$\IS=\IG\big/\bigl(\IG\cap \J\bigr)$.

{}From the commutation relations
defining the algebra $\S$ it follows that $\IS$ is a {\it central}
subalgebra of $\S$. In particular, it is commutative.
The algebra $\IS$ is a natural domain for the Weil homomorphism,
in the regular case.

\subsection{Some Variations}
The construction described in the previous subsection can be
incorporated into a more general context, including some
situations when the (calculus on the) bundle does not admit regular
connections.

We shall first introduce a sufficient condition for the calculus, which
allows that the construction of the Weil homomorphism
works essentially in the same way as in the regular case.

Let us consider a quantum principal bundle $P$, endowed with a
differential calculus $\Omega(P)$. For each $\omega\in\con(P)$ let
$\cal{K}_\omega\subseteq\Omega(P)$ be an ideal generated by
elements of the form
\begin{equation}
v_{\varphi,\omega}=\langle\omega,\varphi\rangle-(-1)^{\partial\varphi}
\langle\varphi,\omega\rangle,
\end{equation}
where $\varphi$ is an arbitrary pseudotensorial form of the type $\adj$.
Evidently, $\cal{K}_\omega$ is *-invariant.

Let us assume that the embedded differential map $\delta$ is
such that there exists a real affine space
$\spl$ consisting of connections $\omega$ satisfying
\begin{equation}\label{spl}
\hor(P)\cap \cal{K}_\omega=\{0\}.
\end{equation}

Evidently, if the bundle admits regular (and multiplicative)
connections then we can choose $\spl=\r(P)$.
However, condition \eqref{spl} is much weaker then the regularity
assumption. In particular, if the embedded differential $\delta$
vanishes identically then $\cal{K}_\omega=\{0\}$ and \eqref{spl} holds
automatically. This covers quantum line bundles, for example.

\begin{pro} \bla{i} If $\omega\in\spl$ then $W^\omega(\vartheta)$ is
closed, for each $\vartheta\in\IG$.

\smallskip
\bla{ii} The cohomological class of $W^\omega(\vartheta)$ in $\Omega(M)$
is independent of the choice of $\omega\in\spl$, and a map
$W\colon\IG\rightarrow H(M)$ is a *-homomorphism.
\end{pro}
\begin{pf} Although the proof is similar as in regular case \cite{D2},
we include it here for completeness. Applying the structure equation we find
$$ (d\omega)^\otimes(\vartheta)=W^\omega(\vartheta)+\sum_{\alpha\beta\eta}
\alpha\langle\omega,\omega\rangle(\eta)\beta $$
for some $\alpha,\beta\in\Omega(P)$ and $\eta\in\Gamma_{\inv}$.
Differentiating both sides we obtain
\begin{multline*}
dW^\omega(\vartheta)=-\sum_{\alpha\beta\eta} \Bigl\{d\alpha
\langle\omega,\omega\rangle(\eta)\beta+(-1)^{\partial\alpha}\alpha
\langle\omega,\omega\rangle(\eta)d\beta\Bigr\}\\
-\sum_{\alpha\beta\eta}(-1)^{\partial\alpha}\alpha\Bigl[
\langle d\omega,\omega\rangle(\eta)-\langle \omega,d\omega\rangle(\eta)
\Bigr]\beta.
\end{multline*}

If $\omega\in\spl$ then the above equality is possible only if
$dW^\omega(\vartheta)=0$, which proves the first statement. Let us consider
another connection $\tau\in\spl$, and let $\omega_t=\omega+t(\tau
-\omega)$ be the segment in $\spl$ determined by $\omega$ and $\tau$.
We have then
$$ \frac{d}{dt} W^{\omega_t}(\vartheta)=dw(\vartheta)+\sum_{f\!g\eta}
f(\langle \varphi,\omega_t\rangle+\langle\omega_t,\varphi
\rangle)(\eta)g+\sum_{uv\zeta}u(\langle d\omega_t,\omega_t\rangle-
\langle\omega_t,d\omega_t\rangle)(\zeta)v, $$
where $w(\vartheta)\in\Omega(M)$ and $f,g\in\hor(P)$ are expressed in
terms of $R_{\omega_t}$, while $\varphi=\omega-\tau$. The splitting
assumption \eqref{spl} implies that the sum of $\cal{K}_\omega$-terms
vanishes. Integrating the above equality from $0$ to
$1$ we conclude that $W^\omega$ and $W^\tau$ coincide in $\Omega(M)$
up to a coboundary.
\end{pf}

In contrast to the regular case however, the Weil homomorphism is
generally not factorizable to the corresponding symmetric algebra.
Another important difference between the regular and the described case
is that the Weil homomorphism $W$ is not longer $\HZ(M)$-valued.

In the context of general bundles and connections on them, the covariant
derivative is not hermitian, and does not satisfy the graded Leibniz
rule. However, appropriate generalizations of these properties hold in
the general case. A deviation of regularity of an arbitrary
connection $\omega$ is naturally expressed via the operator
$\ell_\omega\colon\Gamma_{\inv}\times\hor(P)
\rightarrow\hor(P)$ defined by
\begin{equation}
\ell_\omega(\vartheta,\varphi)=\omega(\vartheta)\varphi-(-1)^{\partial
\varphi}\sum_k\varphi_k\omega(\vartheta\circ c_k).
\end{equation}
With the help of this operator, it is possible to formulate
``hermicity'' and the ``Leibniz rule'' for an arbitrary covariant
derivative.
\begin{pro} \bla{i} The following identities hold
\begin{align*}
D_{\omega}(\varphi\psi)&=D_{\omega}(\varphi)\psi+(-1)^{\partial\varphi}
\varphi D_\omega(\psi)+(-1)^{\partial\varphi}\sum_k\varphi_k\ell_\omega
(\pi(c_k),\psi)\\
D_{\omega}(\varphi)^*&=D_{\omega}(\varphi^*)+\sum_k\ell_\omega
(\pi\k(c_k)^*,\varphi_k^*),
\end{align*}
where $\varphi,\psi\in\hor(P)$.

\smallskip
\bla{ii}  Let us consider a map $\rho\colon\Gamma_{\inv}\otimes\cal{A}
\rightarrow S_{\inv}^{\wedge 2}$ given by
$$
\rho(\vartheta,a)=\delta(\vartheta)\circ a-\delta(\vartheta\circ a)-
(\vartheta\circ a^{(1)})\otimes\pi(a^{(2)})-\pi(a^{(1)})\otimes
(\vartheta\circ a^{(2)}).
$$
Let $u_\omega\colon\Gamma_{\inv}\otimes\cal{A}\rightarrow\hor(P)$ be a
linear map obtained by composing $\omega^\otimes$ and $\rho$. Then the
following identity holds
\begin{multline}\label{Dlw}
D_\omega\ell_\omega(\vartheta,\varphi)=R_\omega(\vartheta)\varphi-
\sum_k\varphi_k R_\omega(\vartheta\circ
c_k)-\ell_\omega(\vartheta,D_\omega(\varphi))\\
{}+\sum_k\varphi_ku_\omega(\vartheta,c_k)+
\sum_j\ell_\omega(\vartheta_j^1,\ell_\omega(\vartheta_j^2,\varphi)),
\end{multline}
where $\delta(\vartheta)=\Sum_j\vartheta_j^1\otimes\vartheta_j^2$.
\end{pro}
\begin{pf} Here we shall prove property \bla{ii}. The statement \bla{i}
follows by modifying trivially a proof given in \cite{D2}. A direct
computation gives
\begin{equation*}
\begin{split}
\langle\omega,\omega\rangle(\vartheta)\varphi&=\sum_j\omega(\vartheta^1_j)
\ell_\omega(\vartheta^2_j,\varphi)+(-1)^{\partial\varphi}\sum_{kj}\omega
(\vartheta_j^1)\varphi_k\omega(\vartheta_j^2\circ c_k)\\
&=\sum_j\ell_\omega(\vartheta_j^1,\ell_\omega(\vartheta_j^2,\varphi))-
(-1)^{\partial\varphi}\sum_{jk}\ell_\omega(\vartheta_2^j,\varphi_k)\omega
(\vartheta_1^j\circ c_k)\\
&\phantom{=}+(-1)^{\partial\varphi}\sum_{kj}
\ell_\omega(\vartheta^1_j,\varphi_k)\omega(\vartheta^2_j\circ
c_k)+\sum_{kj}\varphi_k\omega(\vartheta^1_j\circ
c_k^{(1)})\omega(\vartheta^2_j\circ c_k^{(2)})\\
&=\sum_j\ell_\omega(\vartheta_j^1,\ell_\omega(\vartheta_j^2,\varphi))
-(-1)^{\partial\varphi}\sum_{kl}\ell_\omega(\vartheta_l,\varphi_k)
\omega\bigl[\pi(d_l)\circ c_k\bigr]\\
&\phantom{=}+\sum_{kj}\varphi_k\omega(\vartheta^1_j\circ
c_k^{(1)})\omega(\vartheta^2_j\circ c_k^{(2)}),
\end{split}
\end{equation*}
where $\Sum_j\vartheta_2^j\otimes\vartheta_1^j=\sigma\delta(\vartheta)$
and $\adj(\vartheta)=\Sum_j\vartheta_j\otimes d_j$.
On the other hand
\begin{equation*}
\begin{split}
D_\omega\ell_\omega(\vartheta,\varphi)&=d\ell_\omega(\vartheta,\varphi)+
(-1)^{\partial\varphi}\sum_{jk}\ell_\omega(\vartheta_j,\varphi_k)
\omega\pi(d_jc_k)\\
&=d\omega(\vartheta)\varphi-\omega(\vartheta)d\varphi-(-1)^{\partial\varphi}
\sum_kd(\varphi_k)\omega(\vartheta\circ
c_k)-\sum_k\varphi_kd\omega(\vartheta\circ c_k)\\
&\phantom{=}+(-1)^{\partial\varphi}\sum_{kj}\omega(\vartheta_j)\varphi_k\omega
\bigl(\pi(d_j)\circ c_k\bigr)+(-1)^{\partial\varphi}\sum_k\omega(\vartheta)
\varphi_k\omega\pi(c_k)\\
&\phantom{=}-\sum_{kj}\varphi_k\omega(\vartheta_j\circ
c_k^{(1)})\omega\bigl(\pi(d_j)\circ
c_k^{(2)}\bigr)-\sum_k\varphi_k\omega(\vartheta\circ
c_k^{(1)})\omega\pi(c_k^{(2)}) \\
&=d\omega(\vartheta)\varphi-\omega(\vartheta)D_\omega(\varphi)-
(-1)^{\partial\varphi}\sum_kD_\omega(\varphi_k)\omega(\vartheta\circ c_k)\\
&\phantom{=}-\sum_k\varphi_k\omega(\vartheta\circ
c_k^{(1)})\omega\pi(c_k^{(2)})-\sum_k\varphi_k\omega\pi(c_k^{(1)})\omega
(\vartheta\circ c_k^{(2)})\\
&\phantom{=}+(-1)^{\partial\varphi}\sum_{kj}\ell_\omega(\vartheta_j,\varphi_k)
\omega\bigl(\pi(d_j)\circ c_k\bigr)-\sum_k\varphi_k
d\omega(\vartheta\circ c_k)\\
&=R_\omega(\vartheta)\varphi+\langle\omega,\omega\rangle(\vartheta)\varphi
-\sum_k\varphi_kR_\omega(\vartheta\circ
c_k)-\ell_\omega(\vartheta,D_\omega(\varphi))\\
&\phantom{=}+(-1)^{\partial\varphi}
\sum_{jk}\ell_\omega(\vartheta_j,\varphi_k)\omega\bigl(
\pi(d_j)\circ c_k)\bigr)
-\sum_k\varphi_k\langle\omega,\omega\rangle(\vartheta\circ c_k)\\
&\phantom{=}-\sum_k\varphi_k\omega(\vartheta\circ
c_k^{(1)})\omega\pi(c_k^{(2)})-\sum_k\varphi_k\omega\pi(c_k^{(1)})\omega
(\vartheta\circ c_k^{(2)}).
\end{split}
\end{equation*}
Combining the above expressions, and performing further elementary
transformations we conclude that \eqref{Dlw} holds.
\end{pf}

In general case, for an arbitrary bundle $P$ and a connection $\omega$
on $P$, it is possible to associate to a given calculus $\Omega(P)$ a
new calculus, by passing to an appropriate factoralgebra,
such that the ``projected'' $\omega$ is regular. Let
$\cal{N}_\omega\subseteq\hor(P)$ be a two-sided ideal generated by
elements of the form $\ell_\omega(\vartheta,\varphi)$ and
$$ n_\omega(\vartheta,\varphi)=
R_\omega(\vartheta)\varphi-\sum_k\varphi_kR_\omega(\vartheta\circ c_k)
+\sum_k\varphi_ku_\omega(\vartheta,c_k),
$$
where $\vartheta\in\Gamma_\inv$ and $\varphi\in\hor(P)$.
\begin{pro} The following properties hold
$$ F^\wedge(\cal{N}_\omega)\subseteq\cal{N}_\omega\otimes\cal{A}
\qquad *(\cal{N}_\omega)=\cal{N}_\omega\qquad
D_\omega(\cal{N}_\omega)\subseteq\cal{N}_\omega. $$
\end{pro}
\begin{pf} A direct calculation shows that
\begin{gather*}
\rho(\vartheta,a)\circ b=\rho(\vartheta,ab)-\rho(\vartheta\circ a,b)\\
\rho(\vartheta,a)^*=\rho(\vartheta^*,\k(a)^*)\\
\adj\rho(\vartheta,a)=\sum_l\rho(\vartheta_l,a^{(2)})\otimes
\k(a^{(1)})d_la^{(3)}).
\end{gather*}
It follows that
\begin{gather*}
n_\omega(\vartheta,\varphi)^*=-\sum_k n_\omega(\vartheta^*\circ
\k(c_k)^*,\varphi_k^*)+\bigl\{\mbox{$\ell_\omega$-terms}\bigr\}\\
F^\wedge n_\omega(\vartheta,\varphi)=
\sum_{kl}n_\omega(\vartheta_l,\varphi_k)\otimes d_lc_k.
\end{gather*}
This, together with elementary transformation properties of
$\ell_\omega$ shows that $\cal{N}_\omega$ is $*,F^\wedge$-invariant.

The $D_\omega$-invariance of $\cal{N}_\omega$ follows from
the fact that $D_\omega$ satisfies graded Leibniz rule
modulo $\ell_\omega$-terms, the formula \eqref{Dlw}, as well as from
the formula
for the square \cite{D2} of $D_\omega$. The space $\cal{N}_\omega$ can
be characterized as the minimal $D_\omega$-invariant ideal in $\hor(P)$
containing all the elements $\ell_\omega(\vartheta,\varphi)$.
\end{pf}

Let us consider the factoralgebra
$$\frak{h}(P)=\hor(P)/\cal{N}_\omega.$$ The
*-structure and the right action $F^\wedge$ are projectable from
$\hor(P)$ to $\frak{h}(P)$. Moreover, the covariant derivative
$D_\omega$ projects to a hermitian $F^\wedge$-covariant antiderivation
$D\colon\frak{h}(P)\rightarrow\frak{h}(P)$. Now, we can apply to
$\bigl\{\,\frak{h},F^\wedge,*,D\bigr\}$
a construction of differential calculus presented in \cite{diff}. In such
a way we obtain a graded-differential *-algebra $\Omega_P$,
representing a new calculus on the bundle, a bicovariant *-calculus
$\Psi$ over $G$, and a regular (and multiplicative) connection
$\tau\colon\Psi_\inv\rightarrow\Omega_P$ such that $D_\tau=D$.
However, if the initial connection $\omega$ is sufficiently
``irregular'' then the
the ideal $\cal{N}_\omega$ will be inappropriately large.

\subsection{The Spectral Sequence}
We pass to the study of a spectral sequence naturally associated to
every quantum principal bundle $P$, endowed with a differential
structure $\Omega(P)$. The sequence is based on a filtration of
$\Omega(P)$, induced by the pull back map $\widehat{F}$. The whole
analysis can be viewed as a variation of a classical theme, presented
in \cite{bot} for example.

For each $k\geq 0$ let $\Omega_k(P)\subseteq\Omega(P)$ be the space
consisting of elements having the ``vertical order'' less or equal
$k$. In other words,
$$\Omega_k(P)=\widehat{F}^{-1}\bigl(\Omega(P)\otimes
\Gamma^\wedge_k\bigr)$$
where $\Gamma_k^\wedge$ consists of forms having degrees not exceeding
$k$. These spaces form a filtration of $\Omega(P)$, compatible with the
graded-differential *-structure, in the sense that
$$\Omega_k(P)^*=\Omega_k(P)\qquad d\Omega_k(P)\subseteq\Omega_{k+1}(P)
\qquad\Omega_k(P)=\sideset{}{^\oplus}\sum_{j\geq 0}\Omega_k^j(P).$$
Let us consider a graded-differential *-algebra
$$\kW(P)=\sideset{}{^\oplus}\sum_{k\geq 0}\Omega_k(P), $$
where the grading is given by numbers $k$, and the differential
*-structure is induced from $\Omega(P)$. Let $\kI\colon\kW(P)\rightarrow
\kW(P)$ be the first-order map induced by the inclusions
$\Omega_k(P)\subseteq\Omega_{k+1}(P)$. By definition, this is a
monomorphism of differential *-algebras. Let $\Gr(P)$ be the
graded-differential *-algebra associated to the introduced filtration.
In other words, we have a short exact sequence
\begin{equation}\label{spec-seq}
0@>>> \kW(P) @>{\kI}>> \kW(P) @>>> \Gr(P) @>>> 0
\end{equation}
of differential *-algebras.

The space $\Omega_k(P)$ is linearly spanned by elements of the form
$$w=\varphi\omega(\vartheta_1)\dots\omega(\vartheta_j)$$ where
$\varphi\in\hor(P)$ and $j\leq k$, while $\omega$ is an arbitrary
connection. The algebra $\Gr(P)$ is invariantly isomorphic to the
algebra $\vh(P)$ of ``vertically-horizontally'' decomposed forms.
Explicitly, the isomorphism $\vh(P)\leftrightarrow\Gr(P)$ is
given by
$$ \bigl(\varphi\otimes(\vartheta_1\dots\vartheta_k)\bigr)
\leftrightarrow\bigl[\varphi\omega(\vartheta_1)\dots\omega(\vartheta_k)
+\Omega_{k-1}(P)\bigr].$$
In what follows it will be assumed that the two algebras are identified,
with the help of the above isomorphism.

The factor-differential $d_{v\!h}\colon\vh(P)\rightarrow\vh(P)$
is given by $$d_{v\!h}(\varphi\otimes\vartheta)=
(-1)^{\partial\varphi}\sum_k\varphi_k\otimes\pi(c_k)\vartheta+
(-1)^{\partial\varphi}\varphi\otimes d(\vartheta)$$
where $\Sum_k\varphi_k\otimes c_k=F^\wedge(\varphi)$. The following
technical lemma will prove useful in the computation of the cohomology
algebra of $\vh(P)$.
\begin{lem}\label{lem:triv}
For a given tensorial form $\lambda\colon\Gamma_{\inv}\rightarrow
\Omega(P)$ let us define linear maps
$S_\lambda,T_\lambda\colon\vh(P)\rightarrow\vh(P)$ by the following
equalities
\begin{align*}
S_\lambda\bigl(\varphi\otimes(\vartheta_1\dots\vartheta_n)\bigr)&=
\sum_k\varphi_k\lambda\pi(c_k)\otimes(\vartheta_1\dots\vartheta_n)\\
T_\lambda\bigl(\varphi\otimes(\vartheta_1\dots\vartheta_n)\bigr)&=
(-1)^{\partial\varphi}
\varphi\sum_{i=1}^n(-1)^{i-1}\lambda(\vartheta_{il})\otimes
\Bigl\{\vartheta_1\dots\vartheta_{i-1}\Bigr\}\circ c_{il}
\Bigl\{\vartheta_{i+1}\dots\vartheta_n\Bigr\}
\end{align*}
where $\vartheta_i\in\Gamma_{\inv}$ and $\adj(\vartheta_i)
=\Sum_l\vartheta_{il}\otimes c_{il}$. We have then
\begin{equation}\label{coh-triv}
S_\lambda=T_\lambda d_{v\!h}+(-1)^{\partial\lambda}d_{v\!h}T_\lambda.
\end{equation}
\end{lem}
\begin{pf} A direct computation gives
\begin{multline*}
T_\lambda d_{v\!h}(\varphi\otimes\vartheta)
=(-1)^{\partial\varphi}\sum_k
T_\lambda\bigl(\varphi_k\otimes\pi(c_k)\vartheta
\bigr)+(-1)^{\partial\varphi}T_\lambda
\bigl(\varphi\otimes d(\vartheta)\bigr)\\
=\sum_k\varphi_k\lambda\pi(c_k)\otimes\vartheta+\sum_{ikl}(-1)^i\varphi_k
\lambda(\vartheta_{il})\otimes\Bigl\{\pi(c_k)\vartheta_1\dots\vartheta_{i-1}
\Bigr\}\circ c_{il}\Bigl\{\vartheta_{i+1}\dots\vartheta_n\Bigr\}\\
\phantom{=}+\sum_{i<j}\sum_l(-1)^{i+j}\varphi\lambda(\vartheta_{il})
\otimes\Bigl\{\vartheta_1\dots\vartheta_{i-1}\Bigr\}\circ c_{il}
\Bigl\{\vartheta_{i+1}\dots d(\vartheta_j)\dots\vartheta_n\Bigr\}\\
\phantom{=}+\sum_{i>j}\sum_l(-1)^{i+j-1}\varphi\lambda(\vartheta_{il})
\otimes\Bigl\{\vartheta_1\dots d(\vartheta_j)\dots\vartheta_{i-1}\Bigr\}
\circ c_{il}\Bigl\{\vartheta_{i+1}\dots\vartheta_n\Bigr\}\\
\phantom{=}-\sum_{il}\varphi\lambda(\vartheta_{il})\otimes
\Bigr\{\vartheta_1\dots\vartheta_{i-1}\Bigr\}\circ c_{il}^{(1)}
\Bigr\{\pi(c_{il}^{(2)})\vartheta_{i+1}\dots\vartheta_n\Bigr\},
\end{multline*}
where $\vartheta=\vartheta_1\dots\vartheta_n$. On the other hand
\begin{multline*}
d_{v\!h}T_\lambda(\varphi\otimes\vartheta)=(-1)^{\partial\lambda+i-1}\varphi
\sum_{il}\lambda(\vartheta_{il})\otimes
d\Bigl[\Bigl\{\vartheta_1\dots\vartheta_{i-1}\Bigr\}\circ
c_{il}\Bigl\{\vartheta_{i+1}\dots\vartheta_n\Bigr\}\Bigr]\\
+(-1)^{\partial\lambda+i-1}\sum_{ikl}\varphi_k\lambda(\vartheta_{il})\otimes
\Bigl[\pi(c_k c_{il}^{(1)})\Bigl\{\vartheta_1\dots\vartheta_{i-1}\Bigr\}\circ
c_{il}^{(2)}\Bigl\{\vartheta_{i+1}\dots\vartheta_n\Bigr\}\Bigr].
\end{multline*}
Combining above expressions and performing elementary further
transformations we conclude that \eqref{coh-triv} holds.
\end{pf}

According to the above lemma, the map $S_\lambda$ is cohomologically
trivial (it maps $d_{v\!h}$-cocycles into coboundaries).

The algebra $\vh(P)$ can be naturally decomposed in the following way
\begin{equation}\label{dec-vh}
\vh(P)=\sideset{}{^\oplus}\sum_{u\in\cal{T}}\hbim{u}
\otimes \Lambda_u,
\end{equation}
where $\Lambda_u=H_u\otimes\Gamma_{\inv}^\wedge$ are (understandable
as) subcomplexes, with the differentials specified by
$$d_u(e_i\otimes \eta)=e_i\otimes d(\eta)+\sum_{j=1}^ne_j\otimes
\pi(u_{ji})\eta.$$
\begin{lem} If $u\in\cal{T}$ is such that $\pi(u)\neq 0$ then the
complex $\Lambda_u$ is acyclic.
\end{lem}
\begin{pf}
The maps $S_\lambda$ trivially act on factors $\Lambda_u$, figuring in
the decomposition \eqref{dec-vh}, and they act via composition on spaces
$\hbim{u}$. Therefore, according to Lemma~\ref{lem:triv},
the cohomology of the complex $\Lambda_u$ will be
trivial whenever there exists $\lambda$ satisfying
$S_\lambda\bigl(\hbim{u}\bigr)\neq\{0\}$.
\end{pf}

Let $\infM$ be a space consisting of horizontal forms $\varphi$
satisfying
$$(\id\otimes\pi)F^\wedge(\varphi)=0.$$ In other words,
$$\infM=\sideset{}{^*}\sum_u\cal{H}^u(P)$$
where the summation is performed over $u\in\cal{T}$ satisfying
$\pi(u)=\{0\}$. Evidently, $\Omega(M)\subseteq\infM$.
\begin{lem}
The space $\infM$ is a graded *-subalgebra of $\hor(P)$, closed under
the action of the differential map. \qed
\end{lem}

The elements of $\infM$ play the role of horizontal forms invariant
under infinitesimal vertical diffeomorphisms.

Let $H_{v\!h}(P)$ be
the cohomology algebra of $\vh(P)$. Having in mind that if $\pi(u)=\{0\}$
then the differential $d_u$ reduces to the multiple of  the
orinary differential in $\Gamma_{\inv}^\wedge$ it follows that
\begin{equation}
H_{v\!h}(P)=\infM\otimes H(\Gamma_{\inv}^\wedge).
\end{equation}

Let $E(P)=\Bigl\{E_r(P)\vert r\in\Bbb{N}\Bigr\}$ be the spectral
sequence associated to the short exact sequence
\eqref{spec-seq}.
The introduced filtration of $\Omega(P)$ induces a filtration of the
*-algebra $H(P)$ of cohomology classes. We have $$H_k(P)= \Sum_j^\oplus
H^j_k(P).$$
Applying general theory \cite{bot}, it
follows that the introduced spectral sequence is convergent, and that
$E_\infty(P)$ coincides with the graded *-algebra associated to the
filtered $H(P)$.

By construction, the first cohomology algebra is given by
$E_1(P)=H_{v\!h}(P)$.
\begin{lem}
The differential $d_1$ is given by
\begin{equation}\label{d1}
d_1(w\otimes\vartheta)=dw\otimes\vartheta.
\end{equation}
\end{lem}
\begin{pf} Let us fix a connection $\omega$, and let $\psi=\Sum_{\varphi
\vartheta}\varphi\otimes[\vartheta]^\wedge$ be a $d_{v\!h}$-cocycle. Here
$\varphi\in\infM$ and $\vartheta=\vartheta_1\otimes\dots
\otimes\vartheta_n$, with $\vartheta_i\in\Gamma_{\inv}$.

Then both forms $w=\Sum_{\varphi\vartheta}\varphi\omega(\vartheta_1)\dots
\omega(\vartheta_n)$ and $dw$ belong to $\Omega_n(P)$. By definition,
$d_1[\psi]$ will be represented by the image of $dw$ in $\vh^n(P)$. We have
\begin{multline*}
dw=\sum_{\varphi\vartheta}(-1)^{\partial\varphi}\varphi\sum_{i=1}^n(-1)^{i-1}
\omega(\vartheta_1)\dots\omega(\vartheta_{i-1})\langle\omega,\omega
\rangle(\vartheta_i)\omega(\vartheta_{i+1})\dots\omega(\vartheta_n)\\
+\sum_{\varphi\vartheta}d\varphi\omega^\otimes(\vartheta)
+(-1)^{\partial\varphi}\varphi\sum_{i=1}^n(-1)^{i-1}
\omega(\vartheta_1)\dots\omega(\vartheta_{i-1})R_\omega
(\vartheta_i)\omega(\vartheta_{i+1})\dots\omega(\vartheta_n),
\end{multline*}
and it follows that
$dw-\Sum_{\varphi\vartheta}d\varphi\omega(\vartheta_1)
\dots\omega(\vartheta_n)$ belongs to $\Omega_{n-1}(P)$. Hence,
$$d_1[\psi]=\sum_{\varphi\vartheta}d\varphi\otimes[\vartheta]$$
which completes the proof.
\end{pf}
In particular,
\begin{equation}\label{E2}
E_2(P)=\hinfM\otimes H(\Gamma_{\inv}^\wedge).
\end{equation}

In a special case when $G$ is ``connected,''
only scalar elements of $\cal{A}$ are anihilated by the
differential map, and we have
$$\Omega(M)=\infM\qquad H(M)=\hinfM.$$
Such a connectedness assumption is equivalent to the statement that only
trivial representation $\1$ is anihilated by $\pi$.

\subsection{The Universal Construction}

In this subsection a general approach to quantum characteristic
classes will be presented. It is based on a cohomological interpretation
of the domain of the Weil homomorphism. We shall consider two levels of
generality--quantum principal bundles with general connections, and
structures admitting regular connections.
As first, a characterization the image of the Weil homomorphism will be given.
\begin{pro}
Let us assume that $\omega$ is a regular and multiplicative connection
on $P$. Then the image of the map $W^\omega$ consists of differential
forms on $M$ which are expressible in terms of $\omega$ and $d\omega$.
\end{pro}
\begin{pf}
If $\omega$ is regular and multiplicative then
horizontal elements expressible via $\omega$ and
$d\omega$ are of the form $w=R_\omega^\otimes(\vartheta)$, where
$\vartheta\in\Gamma_{\inv}^\otimes$ is arbitrary. This follows from the
multiplicativity of the horizontal projection $h_\omega$, and the
definition of the curvature map.
If $w\in\Omega(M)$
then it is possible to assume that $\vartheta\in\IG$, because of the
tensoriality of $R_\omega$.
\end{pf}

In generalazing the formalism to the level of arbitrary bundles, we
shall follow {\it the idea of universality}. Algebraic expressions
generating characteristic classes should be the same for all bundles.

In what follows it will be assumed that the higher-order calculus on $G$ is
described by one of the two ``extreme'' algebras--by the universal envelope
$\Gamma^\wedge$ corresponding to the maximal solution, or the
bicovariant \cite{W2} exterior algebra
$\Gamma^\vee$, which describes the minimal appropriate
higher-order calculus. Both cases can be treated essentially in
the same way, and we shall commonly pass from one algebra to another. We
shall denote by $d^\wedge$ and $d^\vee$ the corresponding differentials.

By definition
$$\Gamma^{\vee n}=\Gamma^{\otimes n}/\ker(A_n),$$
where $A_n\colon \Gamma^{\otimes n}\rightarrow\Gamma^{\otimes n}$ are
corresponding braided antisymmetrizers, given by
\begin{equation*}
A_n=\sum_{\pi\in S_n} (-1)^{\pi}\sigma_\pi
\end{equation*}
where operators $\sigma_\pi$ are constructed by replacing
transpositions figuring in minimal decompositions of $\pi$,
by the corresponding $\sigma$-twists (and $\sigma$ is acting on
$\Gamma^{\otimes 2}$ as a bimodule automorphism). Furthermore,
$$ A_{k+l}=(A_k\otimes A_l)A_{kl} \qquad A_{kl}=\sum_{\pi\in S_{kl}}
(-1)^{\pi}\sigma_{\pi^{-1}},$$
where $S_{kl}\subseteq  S_{k+l}$ consists of
permutations preserving the order
of the first $k$  and  the  last $l$ factors. It is interesting to
observe that maps $A_{kl}$ preserve the ideal $S^\wedge$, and hence
there exist natural factorizations
$A^\wedge_{kl}\colon\Gamma^\wedge\rightarrow\Gamma^\wedge$.

The coproduct map is extendible to a homomorphism
$\phi^\vee\colon\Gamma^\vee\rightarrow\Gamma^\vee\grten
\Gamma^\vee$ of graded-differential *-algebras.
By  universality,  there  exists  the  unique  graded-differential
homomorphism  ${\between}\colon \Gamma^\wedge\rightarrow\Gamma^\vee$
reducing to identity maps on $\cal A$ and $\Gamma$. We have
$$\phi^\vee{\between}=({\between}\otimes{\between})\widehat{\phi}.$$

Moreover, all constructions of this paper can be performed also for
``intermediate'' higher-order calculi, described by arbitrary
higher-order graded-differential *-ideals
$S_\star\subseteq\Gamma^\wedge$ satisfying
$$ \widehat{\phi}(S_\star)\subseteq S_\star\grten \Gamma^\wedge +
\Gamma^\wedge\grten S_\star.$$
The above condition implies $S_\star\subseteq[\ker(A)]^\wedge$, which
explaines the minimality of the braided exterior algebra.

It is important to mention that switching from one higher-order calculus
over $G$ to another may influence drastically the properties of the
corresponding characteristic classes for $M$. We shall return to this
point in Section~8.

Let $\WG$ be the universal differential envelope of the tensor algebra
$\Gamma^\otimes_{\inv}$. The algebra $\WG$ will
be endowed with a grading specified by
requiring that elements from $\Gamma_\inv$ have degree $1$.
By definition, $\WG$ is generated (as a differential
algebra) by the space $\Gamma_\inv$ with the only relation $d(1)=0$. The
*-involution on $\Gamma_\inv$ naturally extends to $\WG$, so that
$d\colon\WG\rightarrow\WG$ is a hermitian map. We have
\begin{equation}\label{acycl}
H(\WG)=\Bbb{C},
\end{equation}
which is a property of all universal differential envelopes (of
algebras).

The universality property of $\WG$ implies that there exists
the unique graded-differential homomorphism
$\fWG\colon\WG\rightarrow\WG\grten\Gamma^\wedge$ satisfying
\begin{equation}
\fWG(\vartheta)=\adj(\vartheta)+1\otimes\vartheta
\end{equation}
for each $\vartheta\in\Gamma_{\inv}$. This map is hermitian,
and satisfies the following equality
\begin{equation}
(\fWG\otimes\id)\fWG=(\id\otimes\widehat{\phi})\fWG.
\end{equation}

It is also possible to introduce a natural right action $\rWG\colon
\WG\rightarrow\WG\otimes\cal{A}$, extending the adjoint action map
$\adj$. We have
$$\rWG=(\id\otimes p_*)\fWG. $$

Let $\WGi\subseteq\WG$ be a graded-differential *-subalgebra consisting
of $\fWG$-invariant elements. In other words $w\in\WGi$ iff
$\fWG(w)=w\otimes 1$.

Let $P=(\cal{B},i,F)$ be a quantum principal $G$-bundle over $M$,
endowed with a calculus $\Omega(P)$. Let $\omega$ be an arbitrary
connection on $P$.
\begin{pro}\label{pro:coh-W}
\bla{i} There exists the unique homomorphism
$\widehat{\omega}\colon\WG\rightarrow\Omega(P)$ of differential
algebras satisfying
$\widehat{\omega}(\vartheta)=\omega(\vartheta)$ for each
$\vartheta\in\Gamma_{\inv}$. This map is hermitian, and the diagram
\begin{equation}
\begin{CD}
\WG @>{\mbox{$\fWG$}}>> \WG\grten\Gamma^\wedge \\
@V{\mbox{$\widehat{\omega}$}}VV @VV{\mbox{$\widehat{\omega}\otimes
\id$}}V\\
\Omega(P) @>>{\mbox{$\widehat{F}$}}> \Omega(P)\grten\Gamma^\wedge
\end{CD}
\end{equation}
is commutative. In particular, it follows that
$\widehat{\omega}(\WGi)\subseteq\Omega(M)$.

\smallskip
\bla{ii} The induced cohomology map $W\colon H(\WGi)\rightarrow H(M)$ is
independent of the choice of $\omega$.
\end{pro}
\begin{pf} Property ({\it i\/}) follows from the definition of $\fWG$,
equation \eqref{con}, and the universality of $\WG$. Let $\tau$
be another connection on $P$, and let $L\colon\WG\rightarrow\Omega(P)$
be a linear map specified by
\begin{equation}\label{L}
\begin{gathered}
L(\vartheta)=0\qquad L(d\vartheta)=\varphi(\vartheta)\qquad L(1)=0\\
L(wu)=L(w)\widehat{\tau}(u)+(-1)^{\partial w}\widehat{\omega}(w)L(u),
\end{gathered}
\end{equation}
where $\varphi=\tau-\omega$. A direct computation gives
\begin{equation}
\widehat{\omega}-\widehat{\tau}=Ld+dL,
\end{equation}
and moreover the diagram
\begin{equation}
\begin{CD}
\WG @>{\mbox{$L$}}>> \Omega(P)\\
@V{\mbox{$\fWG$}}VV @VV{\mbox{$\widehat{F}$}}V\\
\WG\grten\Gamma^\wedge
@>>{\mbox{$L\otimes\id$}}> \Omega(P)\grten\Gamma^\wedge
\end{CD}
\end{equation}
is commutative. Hence, $L(\WGi)\subseteq\Omega(M)$. In particular, it
follows that $\widehat{\omega}$ and $\widehat{\tau}$ induce the same
cohomology map $W\colon H(\WGi)\rightarrow H(M)$.
\end{pf}

The constructed map $W$ is a counerpart of Weil homomorphism, at the
level of general quantum principal bundles. Characteristic classes are
therefore labeled by the elements of $H(\WGi)$.

The algebra $\WG$ possesses various properties characteristic to
differential algebras describing the calculus on quantum principal
bundles (however, here there will be no analogs of the base space and
the bundle, because $\WG^0=\Bbb{C}$). In particular, it is possible to
introduce a natural decomposition
\begin{equation}
\WG\leftrightarrow\hh(\WG)\otimes\Gamma_{\inv}^\wedge=\vh(\WG)\qquad
\varphi\iota(\vartheta)\leftrightarrow\varphi\otimes\vartheta
\end{equation}
where $\hh(\WG)\subseteq\WG$ is a graded $*$-subalgebra describing
``horizontal elements,'' defined by
$$\hh(\WG)=\fWG^{-1}(\WG\otimes\cal{A}).$$
It follows that $\rWG[\hh(\WG)]\subseteq\hh(\WG)\otimes\cal{A}$, in other
words $\hh(\WG)$ is $\rWG$-invariant. The algebra $\WGi$ is the
$\rWG$-fixed point subalgebra of $\hh(\WG)$.

Similarly, $\vh(\WG)$ possesses a natural differential *-algebra
structure, in accordance with the analogy with vertically-horizontally
decomposed forms.

Following the ideas of the previous subsection, it is possible
to associate a natural
spectral sequence $E(\WG)=\Bigl\{E_r(\WG)\vert r\in\Bbb{N}\Bigr\}$ to the
algebra $\WG$, so that
$$ E_1(\WG)=H(\vh(\WG),d_{v\!h})=\WGi\otimes H(G)
\qquad E_2(\WG)=H(\WGi)\otimes H(G),$$
and the sequence converges to the trivial cohomology algebra
$E_\infty(\WG)=H(\WG)=\Bbb{C}$. In particular, if $H(G)=\Bbb{C}$ then $\WGi^+$
will be acyclic, too. In this particular case there will be no intrinsic
characteristic classes (the level of general bundles and differential
structures on them).

The following is a prescription of constructing the cocycles from
$\WGi$. Every cocycle $w\in\WGi^+$ is of the form $w=d\varphi$, where
$\varphi$ is some $\rWG$-invariant element of $\WG$. Then we have the
equivalence
$$ w\in\WGi\iff d\bigl\{\fWG(\varphi)-\varphi\otimes 1\bigr\}=0.$$

For every quantum principal bundle $P$, the cocycles representing
characteristic classes are exact, as classes on the bundle, with
$\widehat{\omega}(w)=d\widehat{\omega}(\varphi)$. Therefore, the above
introduced elements $\varphi$ play the role of {\it universal} Chern-Simons
forms. At the level of cohomology classes,

\begin{lem} Let $\cal{C}\subseteq\WG\otimes\Gamma^{\wedge+}$ be the a
subcomplex spanned by elements of the form
$c=\fWG(\varphi)-\varphi\otimes 1$, where $\varphi\in\WG$ is
$\rWG$-invariant. Then the following natural correspondence holds
\begin{equation}\label{WGi-C}
H^n(\WGi)\leftrightarrow H^{n-1}(\cal{C}).
\end{equation}
\end{lem}
\begin{pf} Let us consider the following natural short exact sequence of
complexes
\begin{equation*}
0 @>>> \WGi @>>> I(\WG) @>>> \cal{C} @>>> 0
\end{equation*}
where $I(\WG)\subseteq\WG$ is the subalgebra of $\rWG$-invariant
elements. Passing to the corresponding long exact sequence, and
observing that $H[I(\WG)]=\Bbb{C}$ we conclude that \eqref{WGi-C}
holds.
\end{pf}

Let us compute explicitly {\it the cocycles} of $\WGi$ in dimensions
$2,3$ and $4$. Clearly, closed elements of
$\WGi^2$ are of the form $w=d\vartheta$, where
$\vartheta\in\Gamma_{\inv}$ is $\adj$-invariant and closed in
$\Gamma_{\inv}^\wedge$. It is interesting to observe that if the
higher-order calculus on $G$ is described by the braided exterior
algebra $\Gamma^\vee$, the elements of $H^2(\WGi)$ are labeled by all
$\adj$-invariant elements of $\Gamma_{\inv}$, since such elements are
automatically closed.

Furthermore, closed elements of $\WGi^3$
are generated by Chern-Simons forms
$$\varphi=\Sum_{\vartheta\eta}\vartheta\otimes\eta,$$
where $\vartheta,\eta\in\Gamma_{\inv}$, such that
$$ d^\wedge[\varphi]^\wedge=0\qquad\sigma(\varphi)=\varphi.$$

In general, the corresponding cocycles will be non-trivial
in $\WGi$, which is
an interesting purely quantum phenomena. There exist
characteristic classes in {\it odd} dimensions. However,
if the higher-order calculus on $G$ is described by
$\Gamma^\vee$ then $H^3(\WGi)=\{0\}$.

Third-order Chern-Simons forms can be written as
$$ \varphi=\psi+\sum_{\vartheta\eta} R(\vartheta)\eta,$$
with $\adj$-invariant components
$\Sum_{\vartheta\eta}\vartheta\otimes\eta
\in\Gamma_{\inv}^{\otimes2}$ and $\psi\in\Gamma_{\inv}^{\otimes3}$.
In the above formula
\begin{equation}\label{R}
R(\vartheta)=d\vartheta-\delta(\vartheta)
\end{equation}
is the {\it universal curvature form}. A straightforward
calculation shows that requiring the horizontality of
$w=d\varphi$ is equivalent to the following system of equations
\begin{gather*}
\sum\delta(\vartheta)\otimes\eta=\sum xy\otimes z\qquad\quad
d^\wedge[\psi]^\wedge=0\\
\sum \Bigl\{c^\top(\vartheta)\eta+\vartheta\otimes d^\wedge(\eta)\Bigr\}+
\sum p\otimes qr=0\\
\sum c^\top(p)qr+p\otimes d^\wedge(qr)=0\\
\sum \Bigl\{c^\top(xy)z+xy\otimes d^\wedge(z)\Bigr\}+\sum
\delta(p)\otimes qr=0,
\end{gather*}
where $\adj(\vartheta)=\Sum_k\vartheta_k\otimes c_k$ and
$$ \sum x\otimes y\otimes z=A_{21}(\psi)\qquad\sum p\otimes q\otimes r
=A_{12}(\psi). $$

Now we give the complete description of the horizontal algebra
$\hh(\WG)$. The algebra $\hh(\WG)$ is closed under the operation
\begin{equation}
\ell_\star(\vartheta,\varphi)=\vartheta\varphi-(-1)^{\partial
\varphi}\sum_k\varphi_k(\vartheta\circ c_k),
\end{equation}
where $\vartheta\in\Gamma_{\inv}$ and $\Sum_k\varphi_k\otimes c_k=
\rWG(\varphi)$.

\begin{pro} \bla{i} Let us assume that the higher-order calculus on $G$
is described by the universal envelope $\Gamma^\wedge$. Then $\hh(\WG)$
is the minimal $\ell_\star$-invariant subalgebra of $\WG$ containing
$R(\vartheta)$ and the elements from $S_{\inv}^{\wedge 2}$.

\smallskip
\bla{ii} If $\Gamma^\vee$ describes the higher-order calculus on $G$
then $\hh(\WG)$ is linearly generated by ``elementary'' horizontal forms
$w\in\Gamma^\otimes_{\inv}$ characterized by
\begin{equation}
A_{*1}(w)=0,
\end{equation}
together with the elements obtained from $R(\vartheta)$ by the multiple
actions of $\ell_\star$.
\end{pro}
\begin{pf} Essentially the same reasoning as in \cite{D2}-Section~3,
where the structure of horizontal forms was analyzed in details.
\end{pf}

Let $\sqcup\colon\WG\rightarrow\WG\otimes\Gamma_{\inv}$ be the ``vertical
contraction'' map, obtained by naturally projecting the values of $\fWG$ on
$\WG\otimes\Gamma_{\inv}$. Then we have
\begin{equation}
\hh(\WG)=\ker(\sqcup),
\end{equation}
under the assumption that the whole calculus on $G$ is described by
$\Gamma^\vee$. This means that the system of
equations describing Chern-Simons forms reduces to the first-order
equation (on the second degree).

In any case, these equations always appear as algebraic expressions in
spaces $\Gamma_{\inv}^{\otimes}\otimes\Gamma_{\inv}^{\wedge,\vee}$,
involving maps $\sigma,\delta,\adj,c^\top$ and the differential
$d^{{\wedge}{\vee}}$ in the second factor.

Let us return to the regular case.
Essentially the same cohomological description of the domain of the
Weil homomorphism is possible for bundles admitting regular (and
multiplicative) connections. The only difference is that the algebra
$\WG$ should be factorized through the appropriate ideal, which takes
into account the characteristic property of regular (and multiplicative)
connections.

Let $\WGJ\subseteq\WG$ be the ideal generated by elements from
$S_{\inv}^{\wedge2}$, and the elements of the form
\begin{gather}
j_3(\eta,\vartheta)=\ell_\star(\eta,R(\vartheta))=\eta
R(\vartheta)-\sum_k R(\vartheta_k)(\eta\circ c_k)\\
j_4(\eta,\vartheta)=R(\eta)R(\vartheta)-\sum_k R(\vartheta_k)R
(\eta\circ c_k),
\end{gather}
where $\eta,\vartheta\in\Gamma_{\inv}$ and
$\adj(\vartheta)=\Sum_k\vartheta_k\otimes c_k$.

It follows that $\Gamma_{\inv}^\wedge$ is a subalgebra of
$\WGr=\WG/\WGJ$, in a natural manner. We have
\begin{gather}
\fWG R(\vartheta)=\sum_kR(\vartheta_k)\otimes c_k\\
R(\vartheta)^*=R(\vartheta^*).
\end{gather}
It is easy to see that the ideal $\WGJ$ is $\fWG,*$-invariant. Moreover,
\begin{lem} The ideal $\WGJ$ is $d$-invariant.
\end{lem}
\begin{pf} We have to check that the relations defining $\WGJ$ are
compatible with the action of $d$. Let us first observe that
\begin{equation}\label{dRJ}
dR(\vartheta)=\sum_k R(\vartheta_k)\pi(c_k)+w
\end{equation}
where $w\in\WGJ$. Also, the first argument in $j_3$
can be extended to the whole $\Gamma_{\inv}^\wedge$. A direct
computation gives
\begin{equation*}
\begin{split}
dj_3(\eta,\vartheta)&=w_1+d\eta R(\vartheta)-\sum_k
R(\vartheta_k)d(\eta\circ c_k)
-\sum_k\eta R(\vartheta_k)\pi(c_k)\\
&\phantom{=}-\sum_kR(\vartheta_k)\pi(c_k^{(1)})(\eta\circ
c_k^{(2)})\\
&=w_1+j_4(\eta,\vartheta)+d^\wedge\eta R(\vartheta)-\sum_k
R(\vartheta_k)d^\wedge(\eta\circ c_k)
-\sum_k\eta
R(\vartheta_k)\pi(c_k)\\
&\phantom{=}-\sum_kR(\vartheta_k)\pi(c_k^{(1)})(\eta\circ
c_k^{(2)})\\
&=w_1+j_4(\eta,\vartheta)+j_3(d^\wedge\eta,\vartheta)+\sum_k
\Bigl(R(\vartheta_k)(\eta\circ c_k^{(1)})\pi(c_k^{(2)})-\eta
R(\vartheta_k)\pi(c_k)\Bigr)\\
&=w_1+j_4(\eta,\vartheta)+j_3(d^\wedge\eta,\vartheta)-
\sum_k j_3(\eta,\vartheta_k)\pi(c_k).
\end{split}
\end{equation*}
Applying similar transformations we obtain
\begin{equation*}
\begin{split}
dj_4(\eta,\vartheta)&=w_2+\sum_l R(\eta_l)\pi(d_l)R(\vartheta)
-\sum_k R(\vartheta_k)\pi(c_k^{(1)})R(\eta\circ c_k^{(2)})\\
&\phantom{=}+R(\eta)\sum_k R(\vartheta_k)\pi(c_k)
-\sum_{kl} R(\vartheta_k)R(\eta_l\circ c_k^{(2)})\pi
\bigl[\k(c_k^{(1)})d_lc_k^{(3)}\bigr]\\
&=w_3+\sum_{kl}j_4(\eta_l,\vartheta_k)\pi(d_lc_k),
\end{split}
\end{equation*}
where $w_i\in\WGJ$ and $\adj(\eta)=\Sum_l\eta_l\otimes d_l$. Finally,
$$ d\bigl(\pi(a^{(1)})\otimes\pi(a^{(2)})\bigr)=-
j_3\bigl(\pi(a^{(1)}),\pi(a^{(2)})\bigr) $$
for each $a\in\cal{R}$, in other words $d(S_{\inv}^{\wedge})
\subseteq\WGJ$.
\end{pf}

The maps $\fWG,d,*,\rWG$ are hence projectable to $\WGr$. Obviously,
projected maps (that will be denoted by the same symbols)
are in the same algebraic relations as the original ones.
Let us introduce the horizontal part
$\hh(\WGr)=\fWG^{-1}(\WGr\otimes\cal{A})$ of $\WGr$,
which is a $\rWG$-invariant *-subalgebra of $\WGr$. Applying a similar
reasoning as in \cite{D2} it follows that
\begin{lem} \bla{i} The product map in $\WGr$ induces a graded vector
space isomorphism
\begin{equation}
\hh(\WGr)\otimes\Gamma_{\inv}^\wedge\leftrightarrow\WGr.
\end{equation}

\bla{ii} The map $R\colon\Gamma_{\inv}\rightarrow\WGr$ can be uniquely extended
to a unital *-homomorphism $R\colon\S\rightarrow\WGr$. The extended
$R$ maps isomorphically $\S$ onto $\hh(\WGr)$. Moreover,
the diagram
\begin{equation}
\begin{CD}
\S@>{\mbox{$R$}}>> \hh(\WGr)\\
@V{\mbox{$\adj_{\S}$}}VV @VV{\mbox{$\rWG$}}V\\
\S\otimes\cal{A} @>>{\mbox{$R\otimes\id$}}>
\hh(\WGr)\otimes\cal{A}
\end{CD}
\end{equation}
is commutative. \qed
\end{lem}

Inductively applying \eqref{dRJ} and relations given by $j_3$, it
follows that
$$ d\varphi=\sum_k\varphi_k\pi(c_k) $$
for each $\varphi\in\hh(\WGr)$, where
$\rWG(\varphi)=\Sum_k\varphi_k\otimes c_k$. Let $\WGri\subseteq\WGr$ be
the $\fWG$-fixed-point subalgebra (equivalently, the $\rWG$-fixed-point
subalgebra of $\hh(\WGr)$). From the above formula it follows that
$$d(\WGri)=\{0\},$$
and hence
$$ H(\WGri)=\WGri=\IS $$
in a natural manner. This gives a connection with the previous
definition of the Weil homomorphism.
Essentially the same reasoning as in the proof of
Proposition~\ref{pro:coh-W} leads to the following result.

\begin{pro}
Let us consider a quantum principal $G$-bundle $P=(\cal{B},i,F)$
over $M$, with a calculus $\Omega(P)$ admitting regular and
multiplicative connections. Let $\omega\colon\Gamma_{\inv}
\rightarrow\Omega(P)$ be an arbitrary connection of this type.

\smallskip
\bla{i} There exists the unique homomorphism
$\widehat{\omega}\colon\WGr\rightarrow\Omega(P)$ of differential
algebras extending the identity map on $\Gamma_{\inv}$.
The map $\widehat{\omega}$ is hermitian, and the diagram
\begin{equation}
\begin{CD}
\WGr @>{\mbox{$\fWG$}}>> \WGr\grten\Gamma^\wedge \\
@V{\mbox{$\widehat{\omega}$}}VV @VV{\mbox{$\widehat{\omega}\otimes
\id$}}V\\
\Omega(P) @>>{\mbox{$\widehat{F}$}}> \Omega(P)\grten\Gamma^\wedge
\end{CD}
\end{equation}
is commutative. In particular,
$\widehat{\omega}(\WGri)\subseteq\Omega(M)$.

\smallskip
\bla{ii} The following identities hold
\begin{align}
\widehat{\omega}R&=R_\omega\\
\widehat{\omega}D&=D_\omega\widehat{\omega}\label{wD=Dww}
\end{align}
where $D\colon\WGr\rightarrow\WGr$ is a first-order antiderivation specified
by $$DR(\vartheta)=0\qquad D\vartheta=R(\vartheta).$$

\bla{iii} The induced cohomology map $W\colon H(\WGri)\rightarrow H(M)$ is
independent of the choice of $\omega$. \qed
\end{pro}

The above introduced map $D$ is the {\it
universal covariant derivative},
in accordance with \eqref{wD=Dww}. We have $D^2=0$.

Let $d_{v\!h}\colon\WGr\rightarrow\WGr$ be ``the universal''
vertical differential. By definition, this map is acting in the
following way
\begin{equation}
d_{v\!h}(\varphi\otimes\vartheta)=\sum_k\varphi_k\otimes
\pi(c_k)\vartheta+\varphi\otimes d^\wedge(\vartheta),
\end{equation}
where $\rWG(\varphi)=\Sum_k\varphi_k\otimes c_k$.
\begin{lem}\label{lem:3d} The following identities hold
\begin{align}
Dd_{v\!h}+d_{v\!h}D&=0\\
D+d_{v\!h}&=d.
\end{align}
\end{lem}
\begin{pf} It is sufficient to check the first identity on
elements of the form $\vartheta, R(\vartheta)$ where
$\vartheta\in\Gamma_{\inv}$. We have
\begin{multline*}
d_{v\!h}D(\vartheta)=\sum_kR(\vartheta_k)\otimes\pi(c_k)=
(R\otimes\id)\sigma\delta(\vartheta)-(R\otimes\id)\delta(\vartheta)\\
=(\id\otimes R)\delta(\vartheta)-(R\otimes\id)\delta(\vartheta)=
-Dd_{v\!h}(\vartheta).
\end{multline*}
Furthermore
$$ Dd_{v\!h}R(\vartheta)=\sum_kD\bigl(R(\vartheta_k)\otimes\pi(c_k)
\bigr)=\sum_kR(\vartheta_k)R\pi(c_k)=0=-d_{v\!h}DR(\vartheta). $$
Finally, the second equality restricted to
$\Gamma_{\inv}$ gives us the definition of $R$.
\end{pf}

Let us observe that the map $D$ acts skew-diagonally, with respect to a
natural bigrading in $\WGr$. This implies that
$$H_{\!D}(\WGr)=\sideset{}{^\oplus}\sum_{k,l\in\Bbb{Z}}H^{kl}(\WGr)$$
with the induced bigrading.

To conclude this section, let us analyze the relation between cohomology
algebras $H(\WGi)$ and $H(\WGri)=\WGri$. Let
$p_{\WGJ}\colon\WG\rightarrow\WGr=\WG/\WGJ$ be the factor
projection map.
This map intertwines the corresponding $\fWG$. In particular,
$p_{\WGJ}(\WGi)\subseteq\WGri$.

In the general case two cohomology algebras
will be related through the long exact sequence
\begin{equation}\label{JgrJ}
@>>> H(\WGJ\cap\WGi) @>>> H(\WGi) @>>> H(\WGri) @>>>
H(\WGJ\cap \WGi) @>>>
\end{equation}
associated to
\begin{equation*}
0 @>>> \WGJ\cap\WGi @>>> \WGi @>>> \WGri @>>> 0.
\end{equation*}
The ideal $\WGJ$ can be decomposed as
\begin{equation}
\WGJ\leftrightarrow \hh(\WGJ)
\otimes\Gamma_{\inv}^\wedge\qquad
\hh(\WGJ)=\WGJ\cap\hh(\WG).
\end{equation}

Let us consider the spectral sequence $E_{\WGJ}$
associated to $\WGJ$. From the
form of the second cohomology algebra $E^2_{\WGJ}
=H(\WGJ\cap\WGi)\otimes H(G)$ we conclude that
\begin{pro}
The first non-trivial components of $H(\WGJ)$ and
$H(\WGJ\cap\WGi)$ coincide. \qed
\end{pro}

As we have seen, $H^{2k-1}(\WGri)=\{0\}$ for each $k\in\Bbb{N}$. On the
other hand $\WGi$ will generally contain non-trivial odd-dimensional
classes. From the exact sequence \eqref{JgrJ} it follows that
\begin{equation}
H^{2k-1}(\WGJ\cap \WGi)\twoheadrightarrow H^{2k-1}(\WGi)\qquad
H^{2k}(\WGJ\cap \WGi)\hookrightarrow H^{2k}(\WGi).
\end{equation}
A particularly interesting special case is when the
universal characteristic classes for the regular and the general case
coincide.

\begin{pro}\label{pro:g=r} The following conditions are equivalent:

\smallskip
\bla{i} The complex $\WGJ\cap\WGi$ is acyclic.\par\smallskip
\bla{ii} The complex $\WGJ$ is acyclic.\par\smallskip
\bla{iii} We have $H(\WGr)=\Bbb{C}$.\par\smallskip
\bla{iv} The projection map $p_{\WGJ}\colon
\WGi\rightarrow\WGri$ induces the isomorphism of cohomology
algebras $H(\WGi)$ and $H(\WGri)$.
\end{pro}
\begin{pf} The previous lemma implies that \bla{i} and \bla{ii} are
equivalent. Considering the short exact sequence
\begin{equation*}
0@>>> \WGJ @>>> \WG @>>> \WGr @>>> 0,
\end{equation*}
and remembering that $H^k(\WG)=\{0\}$ for $k\in\Bbb{N}$ we conclude that
$H^k(\WGr)=H^{k+1}(\WGJ)$ and hence the equivalence between \bla{ii} and
\bla{iii} holds. Finally, looking at the exact sequence \eqref{JgrJ} we
conclude that \bla{iv} and \bla{i} are equivalent.
\end{pf}

The isomorphism between cohomology algebras of $\WGi$ and $\WGri$ allows
the essential technical simplification in constructing quantum
characteristic classes--it is sufficient to perform algebraic
constructions assuming the regular case, dealing with the braided
symmetric polynomials.

The spectral sequence $E(\WG)$ converges to the trivial cohomology
$H(\WG)=\Bbb{C}$. Therefore the spaces $H^k(\WGi)$ are naturally
filtered. The filtration is induced from the sequences of
the form
\begin{equation*}
0 @>>> [\,]^{-1}_j\Bigl\{\im(d_{j-1})\Bigr\} @>>>
H(\WGi) @>>> E_j^{*0}(\WG) @>>> 0
\end{equation*}

The trivial convergence information is insufficient to compute the
cohomology of $H(\WGi)$. However, in various interesting special cases
the spectral sequence degenerates (as in classical geometry), and the
triviality property is sufficient to determine all cohomology classes.

\section{Quantum Chern Character}
Through this section we shall deal with structures admitting regular
and multiplicative connections. The aim is to construct the analogue of
the Chern character. This requires introducing the appropriate $K$-rings.
We shall start from the algebra of horizontal forms, and extract vector
bundles as intertwining bimodules (playing the role of
differential forms on $M$ with values in associated vector bundles).
These bimodules generate a ring $K(M,P)$. The Chern character will be
constructed with the help of the covariant derivative (of an arbitrary
regular and multiplicative connection), following classical analogy and
imposing a formula similar to that of \cite{C1}. In particular,
we shall introduce a canonical ``trace'' in the algebras of
corresponding left module homomorphisms. In such a way we obtain
the Chern character as a right module structure $\Ch_M\colon
\HZ(M)\times K(M,P)\rightarrow \HZ(M)$.

It is important to point out a difference between this $\Ch_M$ and
the Chern character introduced in \cite{C1}. The later is defined on the
standard algebraic $K$-group (generated by finite projective right
$\cal{V}$-modules), and has values in the cyclic homology of $M$. In
particular, it is inherently associated to the space $M$, because it
depends exclusively on the algebra $\cal{V}$.
On the other hand the ring $K(M,P)$ explicitly depends on the bundle
$P$, and also on differential structures on $P$ and $G$ (the later
influences the calculus on $M$).
Furthermore, in contrast to the classical case, $\Ch_M$
is not interpretable as a homomorphism between $K(M,P)$ and $\HZ(M)$,
because generally elements from these two spaces do not commute (as
operators of right/left multiplication acting in
$\HZ(M)$). This leads to a
variety of interesting purely quantum phenomenas. In the last
subsection we sketch a construction of quantum Chern classes.

\subsection{The Canonical Trace}

For each $u\in\RepG$ let $\tr_M\colon\lhom(\hbim{u})
\rightarrow\Omega(M)$ be a linear map defined
via the following diagram
\begin{equation}
\begin{CD} \bim{u}^*\otimes_M\!\bim{u}
@>{\mbox{$\id\otimes A$}}>> \bim{u}^*\otimes_M\!\bim{u}\otimes_M\!
\Omega(M)\\
@A{\mbox{$\map{\btwn^u}$}}AA @VV{\mbox{$\langle\rangle_u^-$}}V\\
\cal{V}\leftrightarrow\bim{\1} @>>{\mbox{$\tr_M(A)$}}>
\cal{V}\otimes_M\Omega(M)\leftrightarrow
\Omega(M)\end{CD}
\end{equation}
In the above diagram, $\tr_M(A)$ is a left $\cal{V}$-module
homomorphism, and hence it can be equivalently understood as an element
of $\Omega(M)$. If $A\colon\bim{u}\rightarrow\hbim{u}$ is a
$\cal{V}$-bimodule homomorphism then
$$ f\tr_M(A)=\tr_M(A)f$$
for each $f\in\cal{V}$. Moreover,
$$\tr_M\bigl[\ghom(\hbim{u})\bigr]\subseteq Z(M).$$

It is worth noticing that
$$\tr_M(\map{f})=\tr(C_u^{-1}f)$$
for each $f\in\Mor(u,u)$.

The following equalities hold
\begin{gather}
\tr_M(A\oplus B)=\tr_M(A)+\tr_M(B)\label{tr1}\\
\tr_M(wA)=\tr_M(A)w,\label{tr3}
\end{gather}
where $A,B\in\lhom(\hbim{u,v})$ and $w\in\Omega(M)$ is understood as the
operator of right multiplication. If moreover
$A\in\hom(\bim{u},\hbim{u})$ then
\begin{equation}
\tr_M(A\otimes B)=\tr_M\bigl(\tr_M(A)B\bigr)\label{tr2}.
\end{equation}
The above equalities directly follow from the definition of $\tr_M$.

\begin{lem}
The pairings $\langle\rangle_u^+\colon\bim{u}\otimes_M\!\bim{u}^*
\rightarrow\cal{V}$ and $\langle\rangle_u^-\colon\bim{u}^*
\otimes_M\!\bim{u}\rightarrow\cal{V}$ induce, in
a natural manner, bimodule isomorphisms
of the form
\begin{equation}\label{lrhom}
\lhom\bigl(\bim{u},\cal{V}\bigr)\leftrightarrow\bim{u}^*\qquad
\rhom\bigl(\bim{u},\cal{V}\bigr)\leftrightarrow\bim{u}^*.
\end{equation}
\end{lem}
\begin{pf} We shall check the first isomorphism, the second one follows
by a similar reasoning. As first, the map $\langle\rangle_u^+$ naturally
induces a bimodule homomorphism $\chi_u^+\colon\bim{u}^*\rightarrow\lhom(
\bim{u},\cal{V})$. This map is injective, because of the identity
$$\sum_k \nu_k \langle\mu_k,\psi\rangle_u^+=\psi$$
where $\Sum_k \nu_k\otimes\mu_k=\map{\btwn^u}(1)$. If $\varphi\colon
\bim{u}\rightarrow\cal{V}$ is an arbitrary left
$\cal{V}$-linear map then the element $\psi=\Sum_k\nu_k\varphi(\mu_k)$
satisfies
\begin{multline*}
\chi_u^+(\psi)(\xi)=\langle\xi,\psi\rangle_u^+=\sum_k
\langle\xi,\nu_k\varphi(\mu_k)\rangle_u^+=
\sum_k\langle\xi,\nu_k\rangle_u^+\varphi(\mu_k)\\
=\sum_k\varphi\Bigl[\langle\xi,\nu_k\rangle_u^+\mu_k\Bigr]=\varphi(\xi),
\end{multline*}
for each $\xi\in\bim{u}$. In other words, $\chi_u^+(\psi)=\varphi$, and
hence $\chi_u^+$ is bijective.
\end{pf}

The symmetry between $u$ and $\bar{u}$ implies that also
\begin{equation}\label{lrhom*}
\lhom\bigl(\bim{u}^*,\cal{V}\bigr)\leftrightarrow\bim{u}\qquad
\rhom\bigl(\bim{u}^*,\cal{V}\bigr)\leftrightarrow\bim{u},
\end{equation}
where the isomorphisms are induced by the same contraction maps. More
generally, for each $u,v\in\RepG$ the following bimodule isomorphisms
hold
\begin{equation}\label{hom-uv}
\bim{u}^*\otimes_M\!\bim{v}\leftrightarrow
\lhom(\bim{u},\bim{v})\qquad\bim{u}\otimes_M\!\bim{v}^*\leftrightarrow
\rhom(\bim{v},\bim{u}).
\end{equation}

Furthermore, the identification
$$
(u\times v)^c=v^c\times u^c\qquad (H_u\otimes H_v)^*=H_v^*\otimes H_u^*
$$
induces the following natural bimodule isomorphism
\begin{equation}\label{dual-prod}
\bim{u\times v}^*\leftrightarrow \bim{v}^*\otimes_M\!\bim{u}^*,
\end{equation}
so that the following identities hold
\begin{align*}
\langle x\otimes y,\psi\otimes\varphi\rangle_{u\times v}^+=&
\langle x\langle y,\psi\rangle_v^+\varphi\rangle_u^+\\
\langle \psi\otimes\varphi,x\otimes y\rangle_{u\times v}^-=&
\langle \psi\langle\varphi,x\rangle_u^-y\rangle_v^-.
\end{align*}

The same analysis is applicable to $\Omega(M)$-bimodules $\hbim{u}$.
If $A\in\lhom(\hbim{u},\hbim{v})$ then its {\it right transposed} operator
$A^\top\!\in\rhom(\hbim{v}^*,\hbim{u}^*)$ is defined in a natural manner.
Similarly, each $A\in \rhom(\hbim{u},\hbim{v})$ induces a map $A_\bot
\!\in\lhom(\hbim{v}^*,\hbim{u}^*)$. Explicitly,
\begin{align}
\langle\varphi,Ax\rangle_v^-&=
\langle A_\bot\varphi,x\rangle_u^-\quad\mbox{for
$A\in\rhom(\hbim{u},\hbim{v})$}\\
\langle Ax,\varphi\rangle_v^+&
=\langle x,A^\top\varphi\rangle_u^+\quad\mbox{for
$A\in\lhom(\hbim{u},\hbim{v})$}.
\end{align}
Maps $\top\colon\lhom(\hbim{u},\hbim{v})\rightarrow
\rhom(\hbim{v}^*,\hbim{u}^*)$ and
$\bot\colon\rhom(\hbim{u},\hbim{v})\rightarrow\lhom(\hbim{v}^*,\hbim{u}^*)$
are bimodule anti-isomorphisms. Moreover, $\top$ and $\perp$ act as mutually
inverse transformations and
\begin{equation}\label{tr-prod}
(AB)_\perp=B_\perp A_\perp\qquad(AB)^\top=B^\top A^\top,
\end{equation}
if $A,B$ are appropriate right/left $\Omega(M)$-module homomorphisms.

In particular, if $A\in\hom(\bim{u},\hbim{v})$ then both transposed
maps are defined.
\begin{defn} An element $A\in\hom(\bim{u},\hbim{v})$ is called {\it
transposable} iff $A_\perp=A^\top$.
\end{defn}

If $A$ is transposable then $A^\top\in\hom(\bim{v}^*,\hbim{u}^*)$ is
transposable too. From \eqref{tr-prod} it follows
that the composition of transposable homomorphisms is again
transposable.

It is important to mention that if $f\in\Mor(u,v)$ then $\map{f}\colon
\hbim{v}\rightarrow\hbim{u}$ is transposable. In this case
\begin{equation}
(\map{f})^\top=\map{(f^\top)}.
\end{equation}

Finally, the tensor product of transposable bimodule homomorphisms
is transposable. This property follows from the identities
\begin{equation}
(A\otimes B)^\top=B^\top\otimes A^\top\qquad(A\otimes B)_\perp
=B_\perp\otimes A_\perp
\end{equation}
where the identification \eqref{dual-prod} is assumed.

\begin{lem} We have
\begin{equation}
\tr_M(AB)=(-1)^{\partial\!A\partial\!B}\tr_M(BA),
\end{equation}
for every $A\in\lhom(\hbim{u})$ and transposable $B\in\ghom(\hbim{u})$.
\end{lem}
\begin{pf} Let us suppose that $A\leftrightarrow\varphi\otimes x$, with
$\varphi\in\bim{u}^*$ and $x\in\hbim{u}$. In terms of this
identification we have
$BA-(-1)^{\partial\!A\partial\!B}AB\leftrightarrow\varphi\otimes B(x)-
B^\top(\varphi)\otimes x$ and hence
$$ \tr_M(BA-(-1)^{\partial\!A\partial\!B}AB)=\langle\varphi
,B(x)\rangle^-_u-\langle B^\top(\varphi),x\rangle^-_u=0.\qed $$
\renewcommand{\qed}{}
\end{pf}

Now, we are going to study analogs of covariant derivative maps
which play an important role in the construction of the Chern character.
We shall denote by $\lder(\hbim{u})$ the space of first-order
linear maps $D\colon\hbim{u}\rightarrow\hbim{u}$ satisfying
\begin{equation}
D(\alpha\psi)=d(\alpha)\psi+(-1)^{\partial\alpha}\alpha D(\psi)
\end{equation}
for each $\psi\in\hbim{u}$ and $\alpha\in\Omega(M)$.
Similarly, $\rder(\hbim{u})$ will refer to first-order
maps $D\colon\hbim{u}\rightarrow\hbim{u}$ with the property
\begin{equation}
D(\psi\alpha)=D(\psi)\alpha+(-1)^{\partial\psi}\psi d(\alpha).
\end{equation}
Finally, we define
\begin{equation}
\der(\hbim{u})=\lder(\hbim{u})\cap\rder(\hbim{u}).
\end{equation}

If $D\in\rder(\hbim{u})$ then the formula
\begin{equation}
d\langle\varphi,\psi\rangle_u^-=
\langle D_\bot(\varphi),\psi\rangle_u^-+(-1)^{\partial\varphi}
\langle\varphi,D(\psi)\rangle_u^-
\end{equation}
consistently and uniquely defines a map $D_\bot\in\lder(\hbim{u}^*)$.
Similarly, for each $D\in\lder(\hbim{u})$ the formula
\begin{equation}
d\langle\psi,\varphi\rangle_u^+=
\langle D(\psi),\varphi\rangle_u^++(-1)^{\partial\psi}
\langle\psi,D^\top(\varphi)\rangle_u^+
\end{equation}
defines a transposed map $D^\top\in\rder(\hbim{u}^*)$. The two
transposed operations are understandable as
mutually inverse, in a natural manner.
In particular, for $D\in\der(\hbim{u})$ both transposed maps are defined.

\begin{defn} We say that $D\in\der(\hbim{u})$ is {\it transposable} iff
$D^\top=D_\perp$.
\end{defn}

If $D$ is transposable then $D^\top\in\der(\hbim{u}^*)$ is transposable
too. Furthermore, if $D_1\in\der(\hbim{u})$ and $D_2\in\der(\hbim{v})$
are transposable then $D\colon\hbim{u\times v}\rightarrow\hbim{u\times
v}$ given by
$$D(\vartheta\otimes\eta)=D_1(\vartheta)\otimes\eta+(-1)^{\partial
\vartheta}\vartheta\otimes D_2(\eta)$$
is also transposable and
(modulo the identification \eqref{dual-prod}) we have
$$D^\top(\varphi\otimes\psi)=D_2^\top(\varphi)\otimes\psi+
(-1)^{\partial\varphi}\varphi\otimes D_1^\top(\psi).$$

If $\omega$ is a regular connection on $P$ then every map
$D_{\omega,u}\colon\hbim{u}\rightarrow\hbim{u}$ is transposable, with
$$ D_{\omega,u}^\top=D_{\omega,\bar{u}}.$$

Let us consider an arbitrary $D\in\lder(\hbim{u})$. The formula
\begin{equation}
\nabla(A)=DA-(-1)^{\partial A}AD
\end{equation}
defines a graded derivation $\nabla$ on $\lhom(\hbim{u})$. If in
addition $D\in\der(\hbim{u})$ then the algebras $\rhom(\hbim{u})$ and
$\ghom(\hbim{u})$ are $\nabla$-invariant.
\begin{lem}\label{lem:dtr}
If $D\in\der(\hbim{u})$ is transposable then
\begin{equation}
d_M\tr_M(A)=\tr_M\nabla(A)
\end{equation}
for each $A\in\lhom(\hbim{u})$.
\end{lem}
\begin{pf} It is sufficient to check the above formula for elements of
the form $A\leftrightarrow\varphi\otimes x$. In this case
$DA-(-1)^{\partial A}AD\leftrightarrow D^{\top}(\varphi)\otimes
x+(-1)^{\partial \varphi}\varphi\otimes D(x)$ and we have
\begin{multline*}
d_M\tr_M(A)=d_M\langle\varphi,x\rangle_u^-=
\langle D^{\top}\varphi,x\rangle_u^-+(-1)^{\partial
\varphi}\langle\varphi,Dx\rangle_u^-\\=\tr_M
\bigl(D^\top(\varphi)\otimes x+(-1)^{\partial
\varphi}\varphi\otimes D(x)\bigr)
=\tr_M(DA-(-1)^{\partial A}AD)=\tr_M\nabla(A).\qed
\end{multline*}
\renewcommand{\qed}{}
\end{pf}

For a fixed $u\in\RepG$, let us consider expressions
\begin{equation}
\Theta_n(A,D)=\tr_M(AD^{2n}),
\end{equation}
where $A\in Z(M)$ and $D\in\der(\hbim{u})$. It is worth noticing that
the map $D^2$ is a bimodule homomorphism. This implies that
$\Theta_n(A,D)\in Z(M)$.
\begin{lem}\label{lem:theta}
\bla{i} Let us assume that $D$ is transposable. We have
\begin{equation}
d[\Theta_n(A,D)]=\Theta_n(d(A),D),
\end{equation}
and hence it is possible to pass to the cohomology classes $\HZ(M)$
in the first argument.

\smallskip
\bla{ii} At the level of cohomology classes, $\Theta_n(A,D)$ does not
depend of the choice of $D$.
\end{lem}
\begin{pf}
The first statement directly follows from Lemma~\ref{lem:dtr}. Let us
consider another transposable $D'\in\der(\hbim{u})$, and let
$$D_t=D+tS$$
be a $1$-parameter family of transposable elements, where $S=D'-D$
is a transposable element of $\ghom(\hbim{u})$. Let us assume that $A$
is closed. We have then
\begin{multline*}
\frac{d}{dt}\Theta_n(A,D_t)=\frac{d}{dt}\tr_M(AD_t^{2n})
=\sum_{i=1}^{2n}\tr_M(AD_t^{i-1}S
D_t^{2n-i})\\=\tr_M\bigl(D_tX-(-1)^{\partial X}XD_t\bigr)=d[\tr_M(X)],
\end{multline*}
where $$X=\sum_{i=1}^n D_t^{2i-2}S D_t^{2n-2i}A.$$
It follows that $\Theta_n(A,D)$ and $\Theta_n(A,D')$
belong to the same cohomology class.
\end{pf}

Let us now fix a regular connection $\omega$ on $P$, and consider covariant
derivatives $D_{\omega,u}\colon\hbim{u}\rightarrow\hbim{u}$.
Let $\Ch^u\colon Z(M)\rightarrow Z(M)$ be a linear map given by
\begin{equation}
\Ch^u=\sum_{k=0}^\infty\frac{1}{(2\pi i)^k k!}\Theta_k(A,D_{\omega,u}).
\end{equation}
Strictly speaking, the above infinite sum should be interpreted
as a formal series in a graded algebra.
\begin{pro}\label{pro:chu} \bla{i} The following identities hold
\begin{gather}
\Ch^{u\oplus v}(A)=\Ch^u(A)+\Ch^v(A)\label{ch1}\\
\Ch^{u\times v}(A\otimes B)=\Ch^v(\Ch^u(A)B).\label{ch2}
\end{gather}
\smallskip
\bla{ii} We have
\begin{equation}\label{D:A}
d[\Ch^u(A)]=\Ch^u (d(A)).
\end{equation}
The induced cohomology maps $\Ch^u\colon \HZ(M)\rightarrow \HZ(M)$ are
independent of the choice of $\omega$.
\end{pro}
\begin{pf}
Property \eqref{ch1} directly follows from \eqref{tr1}. Similarly
\eqref{ch2} follows from \eqref{tr2} and the Leibniz rule for $D^2$. We
have
\begin{equation*}
\begin{split}
\Ch^{u\times v}(A\otimes B)&=\sum_{n=0}^\infty\frac{1}{(2\pi i)^{2n} n!}
\tr_M\bigl((A\otimes B)D_{\omega,u\times v}^{2n}\bigr)\\
&=\sum_{n=0}^\infty\sum_{k=0}^n
\frac{1}{(2\pi i)^{2n} n!}\binom{n}{k}
\tr_M\bigl(AD_{\omega,u}^{2k}\otimes BD_{\omega,v}^{2n-2k}\bigr)\\
&=\sum_{kl=0}^\infty\frac{1}{(2\pi i)^{k+l}} \frac{1}{k!l!}
\tr_M\bigl(\tr_M(AD_{\omega,u}^{2k})B D_{\omega,v}^{2l}\bigr)
=\Ch^v(\Ch^u(A)B).
\end{split}
\end{equation*}
The statement \bla{ii} follows from Lemma~\ref{lem:theta}.
\end{pf}

\subsection{Associated $K$-rings $\&$ Chern Characters}
Let $K(G)$ be the $K$-ring associated to the semiring of
equivalence classes of $\RepG$ (with the induced operations of addition
and multiplication). As we have seen, each representation gives
rise to the cohomology map $\Ch^u\colon
\HZ(M)\rightarrow \HZ(M)$.
Equivalent representations give the same map. Formulas
\eqref{ch1}--\eqref{ch2} imply that maps $\Ch^u$ induce a global map
of the form
\begin{equation}
\Ch\colon \HZ(M)\times K(G)\rightarrow \HZ(M)\qquad
([w],[u])\mapsto[\Ch^u(w)],
\end{equation}
which is a right $K(G)$-module structure on $\HZ(M)$.

The ring $K(G)$ generally contains a lot of redundancy, from the point
of view of the introduced module structure on $\HZ(M)$. We shall
now introduce a relativized analogue of the $K$-ring for $M$.

Let $\Phi(M,P)$ be the parametrized set of bimodules $\hbim{u}$, where
$u\in\RepG$. Let $\sim$ be an equivalence relation on $\Phi(M,P)$ defined
in the following way.

We say that $\hbim{u}\sim\hbim{v}$ iff there exists a transposable
isomorphism $J_{uv}\colon \hbim{u}\rightarrow\hbim{v}$. By definition,
the introduced equivalence relation is compatible with the product and the
sum of bimodules. Let $V(M,P)=\Phi(M,P)/\sim$ be the semiring of
corresponding equivalence classes.

\begin{lem}\label{lem:inv-tr} If $\hbim{u}\sim\hbim{v}$ then
\begin{equation}\label{inv-tr}
\tr_M(A)=\tr_M(J_{uv}AJ_{uv}^{-1})
\end{equation}
for each $A\in\lhom(\hbim{u})$.
\end{lem}
\begin{pf} Modulo the natural identifications $\hbim{u}^*\otimes_M\!
\hbim{u}=\lhom(\hbim{u})$ we have
\begin{equation}\label{Tr=()}
\tr_M(A)=\langle A\rangle_u^-.
\end{equation}
It is sufficient to check the formula \eqref{inv-tr} for operators of
the form $A=\varphi\otimes x$. We have
$$J_{uv}AJ_{uv}^{-1}(y)=\langle
J_{uv}^{-1}(y),\varphi\rangle_u^+J_{uv}(x)=
\langle y,(J_{uv}^{-1})^\top\varphi\rangle_v^+J_{uv}(x)$$
which implies
$$\tr_M(J_{uv}AJ_{uv}^{-1})=\langle
(J_{uv}^{-1})^\top\varphi,J_{uv}(x)\rangle_v^-=
\langle\varphi,x\rangle_u^-=\tr_M(A).\qed$$
\renewcommand{\qed}{}
\end{pf}

Let $K(M,P)$ be the $K$-ring associated to $V(M,P)$.
By construction there exists a natural ring epimorphism
$p_M\colon K(G)\rightarrow K(M,P)$, given by $$ p_M([u])=[\hbim{u}].$$

The construction of the map $\Ch\colon \HZ(M)\times
K(G)\rightarrow \HZ(M)$ can be factorized through $p_M$, in other words
\begin{lem} There exists the unique map $\Ch_M\colon \HZ(M)\times K(M,P)
\rightarrow \HZ(M)$ such that the diagram
\begin{equation}
\begin{CD} \HZ(M)\times K(G) @>{\mbox{$\Ch$}}>> \HZ(M)\\
@V{\mbox{$\id\times p_M$}}VV @VV{\mbox{$\id$}}V\\
\HZ(M)\times K(M,P) @>{\mbox{$\Ch_M$}}>> \HZ(M)\\
\end{CD}
\end{equation}
is commutative.
\end{lem}
\begin{pf}
If $D\in\der(\hbim{u})$ is transposable then $J_{uv}DJ_{uv}^{-1}$ is
also transposable, as follows from transposibility of $J_{uv}$. The
statement now follows from Lemma~\ref{lem:inv-tr}.
In other words, $\HZ(M)$ is a right $K(M,P)$-module,
in a natural manner.
\end{pf}

\begin{defn} The map $\Ch_M$ is called {\it the Chern character} for
$M$, relative to $P$.
\end{defn}

\subsection{Explicit Expressions}\label{subs:R}
We are going to find an explicit formula for the Chern character,
in terms of the curvature map. To simplify the
notation, we shall use the
product symbol $\circ$ for the corresponding $K(M,P)$-module structure
on $\HZ(M)$. For each $A\in\lhom(\hbim{u})$ let $A_u$ be a matrix given
by
\begin{equation}
(A_u)_{ij}=\sum_{lk}(C_u^{-1})_{il}\nu_k(e_l^*)(A\mu_k)(e_j).
\end{equation}
We have then
\begin{equation}\label{trAu}
\tr_M(A)=\tr(A_u)
\end{equation}
where on the right-hand side of the above equality figures the standard
trace, and $\Sum_k\nu_k\otimes\mu_k\leftrightarrow I\in\lhom(\hbim{u})$.
For an arbitrary regular connection $\omega$, let $R_{\omega,u}$ be a
$\hor(P)$-valued matrix defined by
$$ (R_{\omega,u})_{ij}=R_\omega\pi(u_{ij}).$$
\begin{lem} \bla{i} The following identity holds
\begin{equation}\label{ch-R}
[A]\circ\bigl[\hbim{u}\bigr]=\sum_{\alpha}\frac{(-1)^\alpha}{(2\pi
i)^\alpha \alpha!}\bigl[\tr(A_uR_{\omega,u}^{\alpha})\bigr].
\end{equation}
\bla{ii} The algebra $\W(M)$ is invariant under the action of $K(M,P)$.
\end{lem}
\begin{pf} Property \bla{i} follows from the expression
$$ D^2\mu(e_i)=-\sum_j\mu(e_j)(R_{\omega,u})_{ji}$$
and equality \eqref{trAu}. To prove \bla{ii}, we shall find an explicit
formula for the action of $K(M,P)$ on elements of $\W(M)$.
If $\vartheta\in\IG$ then
\begin{equation*}
\begin{split}
[R^\otimes(\vartheta)]\circ[\hbim{u}]&=\sum_\alpha\frac{(-1)^\alpha}{
(2\pi i)^\alpha\alpha!}\Bigl[\sum_{ilkn}
(C_u^{-1})_{li}\nu_k(e_i^*)R^\otimes(\vartheta)\mu_k(e_n)(R^\alpha_{\omega,
u})_{nl}\Bigr]\\
&=\sum_\alpha\frac{(-1)^\alpha}{(2\pi
i)^\alpha\alpha!}\Bigl[\sum_{ilkjn}
(C_u^{-1})_{li}\nu_k(e_i^*)\mu_k(e_j)R^\otimes(\vartheta\circ u_{jn})
(R^\alpha_{\omega,u})_{nl}\Bigr]\\
&=\sum_\alpha\frac{(-1)^\alpha}{(2\pi i)^\alpha\alpha!}\Bigl[\sum_{jln}
(C_u^{-1})_{lj}R^\otimes(\vartheta\circ
u_{jn})(R^\alpha_{\omega,u})_{nl}\Bigr],
\end{split}
\end{equation*}
which evidently belongs to $\W(M)$.
\end{pf}

\subsection{Chern Classes}

In this subsection we present a general construction of quantum Chern
classes. As in the whole section,
we work in the regular picture, however the constructed
universal Chern classes are straightforwardly incorporable into the
context of arbitrary bundles, as far as $H(\WGr)=\Bbb{C}$. This follows
from Proposition~\ref{pro:g=r}. However, for bundles admitting regular
connections, Chern classes can be viewed as invariants of associated vector
bundles, as in classical geometry.

Let $\AG\subseteq\cal{A}$ be the subalgebra of $\ad$-invariant
elements. Let $\Delta\colon\AG\rightarrow \IS_\lambda$ be
a map given by
\begin{equation}\label{chern(a)}
\Delta^\lambda(a)=\exp\Bigl\{\sum_{k\in\Bbb{N}} \lambda^k\frac{
\pi(a^{(1)})\dots\pi(a^{(k)})}{k}\Bigr\},
\end{equation}
where $\IS_\lambda$ is the algebra of $\lambda$-series with
coefficients in $\IS$. We shall assume that $\lambda$ is real.
\begin{lem}
The following identities hold
\begin{align}
\Delta^\lambda(a+b)&=\Delta^\lambda(a)\Delta^\lambda(b)\\
\Delta^\lambda[\k(a)^*]&=[\Delta^{-\lambda}(a)]^*.
\end{align}
\end{lem}
\begin{pf} To prove the first equality, it is sufficient to recall that
$\IS$ is commutative. The second equality follows by a direct
computation:
\begin{multline*}
\Delta^\lambda[\k(a)^*]=
\exp\Bigl\{\sum_{k\in\Bbb{N}} \lambda^k\frac{
\pi[\k(a^{(k)})^*]\dots\pi[\k(a^{(1)})^*]}{k}\Bigr\}\\
=\exp\Bigl\{\sum_{k\in\Bbb{N}} (-\lambda)^k\frac{
\pi(a^{(k)})^*\dots\pi(a^{(1)})^*}{k}\Bigr\}=[\Delta^{-\lambda}(a)]^*.\qed
\end{multline*}
\renewcommand{\qed}{}
\end{pf}

The coefficients in the expansion
\begin{equation}
\Delta^\lambda(a)=\sum_{n\geq 0}\lambda^n c_n(a),
\end{equation}
play the role of actual universal Chern classes.
According to the previous lemma
\begin{align}
c_n(a+b)&=\sum_{k+l=n}c_k(a)c_l(b)\\
c_n(a)^*&=(-1)^nc_n[\k(a)^*].
\end{align}

For each $u\in\RepG$ let us denote by $\chi_u\in\AG$ its
{\it character}. Explicitly,
\begin{equation}
\chi_u=\tr(C_u^{-1}u).
\end{equation}
As in the classical theory, we have
$$
\chi_{u\oplus v}=\chi_u+\chi_v\qquad\chi_{u\times v}=\chi_u\chi_v,
$$
and the characters of irreducible representations form a basis
in the space $\AG$.

Let us assume that $\omega$ is a regular connection on a
quantum principal bundle $P$. Then the maps $c_n$ can be further
factorized through the equivalence relation determining $V(M,P)$.
Indeed, we have
\begin{equation}
W\Delta^\lambda(\chi_u)=\Bigl[\exp\tr_M
\Bigl\{\sum_{k\in\Bbb{N}}(-\lambda)^k
D_{\omega,u}^k\Bigr\}\Bigr],
\end{equation}
and hence the formula
$$ c_n[\hbim{u}]=c_n(\chi_u)$$
consistently defines maps $c_n\colon V(M,P)\rightarrow \HZ(M)$,
playing the role of Chern classes of associated vector bundles.

It is interesting to observe that $\Delta^\lambda(a)$ is
generally infinite as a $\lambda$-series, in contrast to
the classical case.

\section{Quantum Line Bundles}
In this section a structural analysis of quantum line bundles
will be performed, including the study of the internal structure
of general differential calculi on them corresponding to a
$1$-dimensional calculus on the structure group $G=U(1)$.
Finally, we
present a construction of the quantum Euler class for arbitrary quantum
principal bundles admitting regular connections.

If $G=U(1)$ then the algebra $\cal{A}$ is generated by a single unitary
element $u\in\cal{A}$, satisfying $\phi(u)=u\otimes u$.
Irreducible representations of $G$ are labeled by integers
$n\leftrightarrow u^n$.

\subsection{Reconstruction Problematics}

The structural analysis of general quantum principal bundles presented
in \cite{D3} essentially simplifies in the case of quantum line bundles.
Let $P=(\cal{B},i,F)$ be such a bundle over a quantum space $M$.

The decomposition to multiple irreducible submodules is simply
\begin{equation}\label{dec-B}
\cal{B}={\sum_{n\in\Bbb{Z}}}^\oplus\cal{B}_n,\qquad
\cal{B}_n=\Bigl\{\,b\in\cal{B}\vert F(b)=b\otimes u^n\Bigr\}
\end{equation}
so that we have
\begin{equation}
(\cal{B}_n)^*=\cal{B}_{-n}\qquad
\cal{B}_{n+m}=\cal{B}_n\otimes_M\!\cal{B}_m,
\end{equation}
in particular the product map in $\cal{B}$ induces
$\cal{V}$-bimodule isomorphisms of
the form $$\mu_+\colon\cal{E}\otimes_M\!\cal{E}^*\rightarrow\cal{V}\qquad
\mu_-\colon\cal{E}^*\otimes_M\!\cal{E}\rightarrow\cal{V},$$
where $\cal{E}=\cal{B}_1$ and $\cal{E}^*=\cal{B}_{-1}$ correspond to
the line bundle (and its conjugate) canonically associated to $P$.
Both maps are
hermitian, in a natural manner.
The structure of the algebra $\cal{B}$ is encoded in $\Bigl\{\cal{V},
\cal{E},\cal{E}^*,\mu_\pm\Bigr\}$.

In terms of identification \eqref{dec-B} the product in
$\cal{B}$ is the standard tensor multiplication,
factorized through the identifications induced by maps $\mu_\pm$. The
*-structure on $\cal{B}$ is uniquely determined by the *-structure on
$\cal{V}$, and the conjugation between $\cal{E}$ and $\cal{E}^*$.
The associativity of the product in $\cal{B}$ implies that
\begin{equation}\label{mu-mu}
\begin{aligned}
\id\otimes\mu_+&=\mu_-\otimes\id\colon\cal{E}^*\otimes_M\!
\cal{E}^{\phantom{*}}\otimes_M\!
\cal{E}^*\rightarrow\cal{V}\\
\id\otimes\mu_-&=\mu_+\otimes\id\colon\cal{E}^{\phantom{*}}
\otimes_M\!\cal{E}^*\otimes_M\!
\cal{E}^{\phantom{*}}\rightarrow\cal{V}.
\end{aligned}
\end{equation}

Conversely, the above properties are sufficient to
construct the bundle $P$, starting from $\cal{E}$ and maps $\mu_\pm$.
Let us assume that a $\cal{V}$-bimodule $\cal{E}$ is given, and let
$\cal{E}^*$ be {\it the  conjugate} bimodule. Let us assume that
there exist hermitian $\cal{V}$-bimodule isomorphisms
$\mu_-\colon\cal{E}^*\otimes_M\!\cal{E}\rightarrow\cal{V}$ and
$\mu_+\colon\cal{E}\otimes_M\!\cal{E}^*\rightarrow\cal{V}$
such that properties \eqref{mu-mu} hold.

Let us define first bimodules $\cal{B}_n$, and then
the algebra $\cal{B}$ by \eqref{dec-B},
at the level of vector spaces. Let the product in $\cal{B}$ be induced
by the tensor product and identifications given by
maps $\mu_\pm$. Properties \eqref{mu-mu} ensure the associativity of the
introduced product.
By construction, the *-structure on $\cal{V}$ and the mutual conjugation
between $\cal{E}$ and $\cal{E}^*$ admit the unique extension to a
*-structure on $\cal{B}$.

Let $F\colon\cal{B}\rightarrow\cal{B}\otimes\cal{A}$ be a *-homomorphism
given by $F(q)=q\otimes u^n$,
where $q\in\cal{B}_n$. By
construction $\cal{V}$ is the fixed-point subalgebra for this action.
Let $i\colon\cal{V}\hookrightarrow\cal{B}$ be the canonical inclusion
map.
\begin{lem}
The triplet $P=(\cal{B},i,F)$ is a principal $G$-bundle over $M$. \qed
\end{lem}

\subsection{Differential Structures and The Euler Class}

Let us assume that the calculus on $G$ is $1$-dimensional and
non-standard, and that
higher-order components of $\Gamma^\wedge$ are trivial. Calculi $\Gamma$
satisfying the above assumptions are given by ideals of the form
$$ \cal{R}=
\bigl(u^{-1}+\frac{u}{\lambda}-(1+\frac{1}{\lambda})\bigr)\cal{A},$$
where $\lambda\in\Re\setminus\{-1,0,1\}$. The reality of $\lambda$
ensures the existence of the *-structure on $\Gamma$. The element
$\zeta=\pi(u)/(1-\lambda)$ spans the space $\Gamma_{\inv}$, and we have
$$ \zeta\circ a=p_\lambda(a)\zeta$$
where $p_\lambda\colon\cal{A}\rightarrow\Bbb{C}$ is the character
specified by $p_\lambda(u)=\lambda$. The ideal $S_{\inv}^\wedge$ is
generated by $\zeta\otimes\zeta$, and consequently
$\Gamma_{\inv}^\wedge=\Bbb{C}\oplus\Gamma_{\inv}$.

Let $P$ be an arbitrary quantum principal $G$-bundle over $M$ and let
$\Omega(P)$ be a graded-differential *-algebra describing the full
calculus on $P$. Let us fix a connection $\omega$ on $P$.

According to the general theory each element $w\in\Omega(P)$ can be
uniquely decomposed as
$$w=\varphi+\psi\omega$$
where $\varphi,\psi\in\hor(P)$ with $\varphi=h_\omega(w)$, and
we have identified $\omega=\omega(\zeta)$. This enables us to write
\begin{equation}\label{W-split}
\Omega(P)=\hor(P)\oplus\hor(P)
\end{equation}
at the level of left $\hor(P)$-modules.

We pass to the study of the representation of the differential *-algebra
structure on $\Omega(P)$, in terms of the above identification. Let
$\ell\colon\hor(P)\rightarrow\hor(P)$ be a linear map defined by
\begin{equation}
\ell(\varphi)=\omega\varphi-(-1)^{\partial\varphi}\chi(\varphi)
\omega,
\end{equation}
where $\chi\colon\hor(P)\rightarrow\hor(P)$ is an authomorphism given by
\begin{equation}
\chi(\varphi)=(\id\otimes p_\lambda)F^\wedge(\varphi).
\end{equation}
The first-order map $\ell=\ell_\omega(\zeta,*)$ measures the lack
of regularity of $\omega$. It satisfies a $\chi$-relativized
graded Leibniz rule
\begin{equation}\label{rel-Der}
\ell(\varphi\psi)=\ell(\varphi)\psi+(-1)^{\partial\varphi}\chi(\varphi)
\ell(\psi).
\end{equation}
The following covariance properties hold
\begin{equation}
\begin{gathered}
\chi(\varphi)^*=\chi^{-1}(\varphi^*)\\
\ell(\varphi)^*=\ell\chi^{-1}(\varphi^*)
\end{gathered}\qquad
\begin{aligned}
F^\wedge\chi&=(\chi\otimes\id)F^\wedge\\
F^\wedge\ell&=(\ell\otimes\id)F^\wedge.
\end{aligned}
\end{equation}

In terms of the realization \eqref{W-split}, the product in $\Omega(P)$
is given by
\begin{multline}\label{prod}
(w,u)(\varphi,\psi)=\\=\bigl(w\varphi+u\ell(\varphi)
+(-1)^{\partial \psi}u\chi(\psi)\rho,u\ell(\psi)+(-1)^{\partial\varphi}
u\chi(\varphi)+w\psi+(-1)^{\partial\psi}u\chi(\psi)\bigr).
\end{multline}
Here $\rho=\omega^2\in\Omega(M)$ measures the lack of
multiplicativity of $\omega$.
The *-structure on $\Omega(P)$ is specified  by
\begin{equation}\label{*-real}
(w,u)^*=\bigl(w^*-(-1)^{\partial u}\ell(u^*),-\chi(u^*)\bigr).
\end{equation}

Let $D\colon\hor(P)\rightarrow\hor(P)$ be the covariant derivative
associated to $\omega$, and $R=d\omega$ the corresponding
curvature. From the general theory it follows that the following
identities hold
\begin{gather}
D(\varphi\psi)=D(\varphi)\psi+(-1)^{\partial\varphi}\varphi D(\psi)
+(-1)^{\partial\varphi}\bigl(\varphi-\chi(\varphi)\bigr)\ell(\psi)\\
-D^2(\varphi)=(I-\chi)(\varphi)R+(I-\chi)^2(\varphi)\rho.
\end{gather}

Furthermore, the differential map is given by
\begin{multline}\label{d-real}
d(w,u)=\\=\bigl(D(w)+(-1)^{\partial u}uR+(-1)^{\partial u}[u-\chi(u)
]\rho, D(u)+(-1)^{\partial w}[w-\chi(w)]\bigr).
\end{multline}
The pull back map $\widehat{F}$ is representable as follows
\begin{equation}
\widehat{F}(w,u)=\bigl(F^\wedge(w)+F^\wedge(u)\zeta,
F^\wedge(u)\bigr),
\end{equation}
where $F^\wedge\colon\hor(P)\rightarrow\hor(P)\otimes\cal{A}$ is the
right action map.

Finally, the following additional algebraic properties hold:
\begin{equation}
\begin{gathered}
R^*=-R\qquad\rho^*=-\rho\qquad D(\varphi)^*=D(\varphi^*)+
\ell(\varphi^*)-\ell(\varphi)^*\\
\ell(\rho)=0\qquad d\rho=-\ell(R)\qquad
\rho\varphi=\ell^2(\varphi)+\chi^2(\varphi)\rho\\
D\ell(\varphi)+\ell D(\varphi)=R\varphi-\chi(\varphi)R-2\chi\bigl(
\varphi-\chi(\varphi)\bigr)\rho\qquad DR=0.
\end{gathered}
\end{equation}
This completes the structural analysis of $\Omega(P)$. We see that the
calculus is completely determined by maps
$\bigl\{D,*,F^\wedge,\chi,\ell\bigr\}$ acting in $\hor(P)$, and
second-order differential forms $\bigr\{R,\rho\bigr\}$ on the base
manifold.

Conversely, the derived algebraic properties can be taken as a starting
point in constructing the calculus on $P$. Let us assume that a
graded *-algebra $\hor_P$ representing horizontal forms
is given, such that $\hor_P^0=\cal{B}$. Let us further assume that
$\hor_P$ is endowed with maps $D,\chi,\ell\colon\hor_P
\rightarrow\hor_P$ and the extension
$F^\wedge\colon\hor_P\rightarrow\hor_P\otimes\cal{A}$ of the right
action $F$. Let $\Omega_M\subseteq\hor_P$ be the fixed-point
subalgebra for the action $F^\wedge$, and let $R,\rho$ be second-order
elements of $\Omega_M$. Finally, let us assume that all the above
derived algebraic properties hold. Then formulas \eqref{prod}, \eqref{*-real}
and \eqref{d-real} determine a differential *-algebra structure on
$\Omega_P=\hor_P\oplus\hor_P$. Moreover, the formula
\begin{equation}
\omega(\zeta)=(0,1)
\end{equation}
determines a connection $\omega$. The ``reconstruction circle''
is closed by observing that $\bigl\{R,\rho,D,\ell\bigr\}$ can be
re-expressed in terms of $\omega$, in the original way.

Let us analyze the cohomology algebra
$H(\WGi)$ for the given $1$-dimensional $\Gamma$.
By definition, $\WG$ is the free
differential algebra generated by a first order generator $\zeta$. The
map $\fWG$ acts as follows
\begin{equation}
\fWG(\zeta)=\zeta\otimes 1+1\otimes\zeta\qquad\fWG(d\zeta)
=d\zeta\otimes 1.
\end{equation}
The space of $d$-cocycles of $\WG^+$ is spanned by the elements
$w=d(\zeta^{n_1})\dots d(\zeta^{n_k})$.
\begin{lem} We have
\begin{equation}
H^{2k}(\WGi)=\Bbb{C}\qquad H^{2k+1}(\WGi)=\{0\}
\end{equation}
for each $k\in\Bbb{N}\cup\{0\}$. More precisely, the algebra $H(\WGi)$
is generated by the second-order element $d\zeta$.
\end{lem}
\begin{pf}
It follows by analyzing the structure of the differential algebra $E_2(\WG)$,
taking into account that $H(G)=\Bbb{C}\oplus\Bbb{C}$, and observing that
the spectral sequence stabilizes at the third term.
\end{pf}

The algebra of characteristic classes is therefore generated by
the cohomological class $e_M=[d\omega]$ in $\Omega(M)$, which is a
direct counterpart of the Euler class. The connection form $\omega$ is
here a counterpart of the global angular form \cite{bot}.

Let us assume that the bundle admits regular connections $\omega$,
and find the action of $K(M,P)$ on the elements of a subalgebra
$\W(M)\subseteq \HZ(M)$. We have
$$ \pi(u^n)=(1-\lambda^n)\zeta$$
for each $n\in\Bbb{Z}$, and a straightforward application
of the explicit formulas of Subsection~\ref{subs:R} gives
\begin{equation}
(e_M^k)\circ [\hbim{n}]=\lambda^{nk}\sum_{\alpha\geq 0}
\frac{(\lambda^n-1)^\alpha}
{(2\pi i)^\alpha\alpha!}e_M^{k+\alpha}
\end{equation}
where $n$ labels non-trivial irreducible representations.

To conclude this section, we shall present a general construction of
quantum Euler class, at the level of bundles $P$ possesing regular
connections. We follow \cite{KN}, the structure group $G$ is arbitrary.

Let $u\in\RepG$ be a representation satisfying $u\sim u^c$. Let
$S\colon H_u\rightarrow H_u^*$ be an equivalence between
$u$ and $u^{c}$. The formula
$$\lambda_{u,S}(e_i\otimes e_j)=\sum_k S_{ki}u_{kj}\qquad
S(e_i)=\sum_kS_{ki}e_k^*$$
defines a map $\lambda_{u,S}\colon H_u\otimes H_u\rightarrow\cal{A}$
intertwining $u\times u$ and the adjoint action $\ad$.

Let us now assume that there exists a braid operator $\tau\colon
H_u^{\otimes 2}\rightarrow H_u^{\otimes 2}$ which intertwines $u\times
u$, and such that there exists a {\it volume element} in the
corresponding braided exterior algebra $H_u^\wedge$.
In other words, there exists $k\in\Bbb{N}$ such that $H_u^{\wedge
k}$ is $1$-dimensional, while
higher order components vanish. Let us consider
a volume element $w\in H_u^{\wedge k}\setminus\{0\}$. We have
$$ u^{\wedge}(w)=w\otimes\Delta_u$$
where $\Delta_u\in\cal{A}$ is the corresponding quantum determinant and
$u^\wedge$ the corresponding representation of $G$ in $H_u^\wedge$.
It is worth noticing that $k$ is not necessarily equal to the
dimension of $H_u$.

Let us finally assume that $k=2n$ and that
$u$ is {\it unimodular}, in the sense that $\Delta_u=1$.
Let $w_\S\in\IS$ be an element given by
\begin{equation}
w_\S=c_n[\pi^{\otimes n}\lambda_{u,S}^{\otimes n}(w)]_\S,\qquad c_n>0.
\end{equation}
\begin{defn} The characteristic class
\begin{equation}
e_M^u=W(w_\S)
\end{equation}
is called the {\it quantum Euler class} of $P$, relative to $u$.
\end{defn}
\section{Bundles Over Classical Manifolds}

In this section we shall study characteristic classes of locally trivial
bundles over classical smooth manifolds, applying the theory developed
in \cite{D1}. Let $M$ be a compact smooth manifold, $P$ a locally
trivial principal $G$-bundle over $M$, and let us assume that $\Gamma$
is the minimal admissible (bicovariant *-) calculus over $G$. The
admissibility here means that $\Gamma$ is compatible, in appropriate
sense, with all local retrivializations (transition functions)
of the bundle. Let
$\Omega(P)$ be a graded-differential *-algebra canonically associated
\cite{D1} to $P$ and $\Gamma$. The main property of $\Omega(P)$ is its
full {\it local trivializability}, in the sense that all local
trivializations of the bundle locally trivialize the calculus, too.

Let $\tau=\bigl(\pi_U\bigr)_{U\in\cal{U}}$ be an arbitrary
trivialization system for $P$. Locally, connection forms are (uniquely)
expressed via the corresponding ``gauge fields'' as
$$\pi_U^\wedge\omega(\vartheta)=(A_U\otimes\id)\adj(\vartheta)+
1_U\otimes\vartheta,$$
where $A_U\colon\Gamma_{\inv}\rightarrow\Omega(U)$ are hermitian maps,
mutually related in the appropriate way on overlapings of regions from
$\cal{U}$. Here $\pi_U^\wedge\colon\Omega(P)\rightarrow\Omega(U)\grten
\Gamma^\wedge$ are the corresponding trivializations of the calculus.
The curvature is expressed by
\begin{equation}\label{F}
\pi_U^\wedge R_\omega(\vartheta)=(F_U\otimes\id)\adj(\vartheta),
\qquad F_U=dA_U-\langle A_U,A_U\rangle.
\end{equation}
Regular connections are characterized by
\begin{equation}\label{cl-reg}
A_U(\vartheta\circ a)=\e(a)A_U(\vartheta),
\end{equation}
for each $\vartheta\in\Gamma_{\inv}$ and $a\in\cal{A}$.

According to \cite{D1} quantum $G$-bundles $P$ are in a $1$--$1$
correspondence with classical $G_{cl}$-bundles (over the same manifold),
where $G_{cl}$ is {\it the classical part} of $G$ (consisting of
classical points of $G$). Explicitly, the elements of $P_{cl}$ are in a
natural correspondence with classical points of $P$ (the *-characters of
$\cal{B}$).

Regularity condition \eqref{cl-reg} is equivalent to the statement that
$A_U$ are factorizable through the canonical projection map
$\nu\colon\Gamma_{\inv}\rightarrow\lie(G_{cl})^*$. In this
sense regular connections are understandable as connections on $P_{cl}$
(classical connections, in the terminology of \cite{D1}).
Moreover, the curvature operator ``commutes'' with the restriction from
$P$ to $P_{cl}$. This follows from equality
$$F_U(\vartheta\circ a)=\e(a)F_U(\vartheta).$$

In summary, characteristic classes of $P$ are interpretable as
(classical) characteristic classes of $P_{cl}$. The converse is
generally not true, because $\nu$ generally maps invariant elements of
$\S$ into a proper subalgebra of invariant polynomials
over $\lie(G_{cl})$.
The classical nature of quantum characteristic classes is
explained by the correspondence $P\leftrightarrow P_{cl}$. The whole
topological nontriviality of $P$ is already contained in its classical
part $P_{cl}$.

The above described concept of local triviality can be further naturally
generalized in two ways \cite{b,pfl}. One way of generalizing consists
in replacing $M$ with an appropriate quantum space, endowed with a
system of noncommutatitive ``open sets'' $U$, locally trivializing $P$.
Secondly, it is possible to relax the requirement that locally the
bundle and the calculus are tensor products of algebras--it is
sufficient to ask that locally $\Omega(P)\leftrightarrow
\Omega(U)\otimes\Gamma^\wedge$, but with the appropriate
crossed product structure (so that locally $\widehat{F}\leftrightarrow
\id\otimes\widehat{\phi}$). Local expressions for $\omega$ and
$R_\omega$ will be essentially the same in this general context.

Furthermore, requiring that all local trivializations of the bundle
locally trivialize the calculus implies an ``admissibility''
constraint on the calculus $\Gamma$, and there exists the minimal
calculus (for a given bundle $P$) resolving the constraints at the
first-order level. In both choices $\{\Gamma^\wedge,\Gamma^\vee\}$ for
the higher-order calculus on $G$, compatibility with retrivializations
will be automatical, at the higher-order levels.

However, even when the base is classical, the class of locally trivial
bundles is too restrictive. As an interesting example, we shall describe a
locally-nontrivializable quantum line bundle over a classical
compact manifold $M$.

Let us consider an arbitrary quantum line bundle $P$ over $M$, and
assume that $\cal{V}=S(M)$. At the level of left modules, we can write
$\cal{E}=S(\xi)$ in a natural manner, where $\xi$ is a line vector
bundle (in the standard sense) over $M$. In terms of this
identification, the right module structure is described by
$$ \psi f=\varepsilon(f)\psi$$
where $\varepsilon\colon\cal{V}\rightarrow\cal{V}$ is a
*-automorphism. This map further extends to the automorphism $\varepsilon
\colon\Omega(M)\rightarrow\Omega(M)$, where $\Omega(M)$ is the standard
graded-differential *-algebra of smooth differential forms. Actually,
$\varepsilon$ is induced (via the pull back) by a diffeomorphism
$\varepsilon_*\colon M\rightarrow M$.

The algebra $\cal{B}$ is noncommutative, if $\varepsilon\neq I$.
If the line bundle $\xi$ is non-trivial, and if the map
$\varepsilon_*$ is such that $\{\emptyset, M\}$ are the only
$\varepsilon_*$-invariant open subsets of $M$ then the bundle
$P$ is not locally trivial.

This example is also interesting from the point of view of
characteristic classes. Let $\hor(P)$ be the *-algebra of horizontal
forms, built from $\{\xi,\varepsilon,\Omega(M)\}$.

\begin{pro}
There exists a
natural corresponence between regular covariant derivative maps
$D\colon\hor(P)\rightarrow\hor(P)$ and standard connections
$\nabla$ on $\xi$. \qed
\end{pro}

Let $\Phi\in\Omega^2(M)$ be the curvature of $\nabla$. Let $\Gamma$ be
the minimal calculus on $G=U(1)$ compatible with $D$, in the sense of
\cite{diff}. Let $\Omega(P)$ be the corresponding graded-differential
*-algebra representing the calculus on $P$, constructed applying the
construction of \cite{D2}-Section~6. Finally, let $\omega\colon
\Gamma_{\inv}\rightarrow\Omega(P)$ be the canonical regular connection
on $P$, so that we have $D=D_\omega$. The curvature of this connection
is given by
\begin{gather*}
R_\omega\pi(u^n)=\varepsilon^{-1}(\Phi)+\dots+\varepsilon^{-n}(\Phi)\\
-R_\omega\pi(u^{-n})=\Phi+\dots+\varepsilon^{n-1}(\Phi),
\end{gather*}
where $n\in\Bbb{N}$.

In the general case, the resulting calculus $\Gamma$ will be
multi-dimensional. The algebra $\W(M)$ of characteristic classes will be
generated by the above expressions. In other words, it is the minimal
$\varepsilon$-invariant subalgebra of $H(M)$ containing the class of
$\Phi$.

\section{Framed Bundles}

In general, there exist no natural prescriptions of constructing
differential calculi over a given quantum principal bundle $P$. However,
when $P$ is equipped with the appropriate additional structure,
it will be possible to associate a distinguished differential algebra
$\Omega(P)$ to the bundle. A large class of structures of this kind
giving non-trivial calculi is given by framed bundles \cite{frm1,frm2}.
Let us sketch the construction of a natural calculus on such bundles.

Let $\Psi$ be a bicovariant *-bimodule over $G$. According to the
general theory \cite{W2}, it is structuralized as $\Psi\leftrightarrow
\cal{A}\otimes V$, where $V=\Psi_{\inv}$.

Let $v\colon V\rightarrow V\otimes\cal{A}$ and $\circ$ be the natural
right action and the right $\cal{A}$-module structure on $V$,
respectively. Let $\tau$ be the canonical braid operator on
$V\otimes V$. In what follows it will be assumed that $\ker(I+\tau)\neq
\{0\}$.

Let $V^\wedge$ be the corresponding $\tau$-exterior algebra, obtained by
factorizing $V^\otimes$ through the relations $\im(I+\tau)$. The
*-structure, $v$ and $\circ$ are naturally extendible to $V^\wedge$, so
that we have a bicovariant graded *-algebra $\Psi^\wedge\leftrightarrow
\cal{A}\otimes V^\wedge$.

Let us consider a graded vector space
$$\hor_P=\cal{B}\otimes V^\wedge$$
endowed with the following *-algebra structure
\begin{align*}
(q\otimes\vartheta)(b\otimes\eta)&=\sum_iqb_i\otimes(\vartheta\circ c_i)
\eta\\
(b\otimes\vartheta)^*&=\sum_i b_i^*\otimes(\vartheta^*\circ c_i^*),
\end{align*}
where $\Sum_i b_i\otimes c_i=F(b)$.

There exists the unique action $F^\wedge\colon\hor_P\rightarrow\hor_P
\otimes\cal{A}$ extending actions $F$ and $v^\wedge\colon V^\wedge
\rightarrow V^\wedge\otimes\cal{A}$. Let us observe that the above
definition of the product in $\hor_P$, together with the definition of
$V^\wedge$, implies the following commutation relations
$$ \vartheta\varphi=(-1)^{\partial\varphi\partial\vartheta}\varphi_k(
\vartheta\circ c_k),$$
with $\Sum_k \varphi_k\otimes c_k=F^\wedge(\varphi)$ and $\vartheta
\in V^\wedge$.

By definition, the bundle $P$ is $\Psi$-framed, iff there exists a first-order
antiderivation $D\colon\hor_P\rightarrow\hor_P$ satisfying
\begin{gather*}
F^\wedge D=(D\otimes\id)F^\wedge\qquad {*}D=D{*}\\
D(V^\wedge)=\{0\}\qquad D^2(\cal{V})=\{0\},
\end{gather*}
and such that the canonical map $\sharp_D\colon\cal{B}\otimes\cal{V}
\rightarrow\hor^1(P)$ given by $\sharp_D(b\otimes f)=bD(f)$ is
surjective.

Let $\Omega_M\subseteq\hor_P$ be the $F^\wedge$-fixed point subalgebra.
The covariance condition implies $D(\Omega_M)\subseteq\Omega_M$. Let
$d_M\colon\Omega_M\rightarrow\Omega_M$ be the restriction map. By
construction, $\Omega_M^0=\cal{V}$. Furthermore, it follows that
$d_M^2=0$. Moreover it turns out that $\Omega_M$ is generated by
$\cal{V}$, as a differential algebra equipped with $d_M$.

Applying to $\bigl\{\hor_P,F^\wedge,D,\Omega_M\bigr\}$
a general construction presented in \cite{diff} we
obtain in the intrinsic way a graded-differential *-algebra
$\Omega(P)$ representing the calculus on $P$, a bicovariant *-calculus
$\Gamma$ over $G$, and a regular and multiplicative
connection $\omega\colon\Gamma_{\inv}
\rightarrow\Omega(P)$ such that $\hor_P$ is recovered as the
corresponding *-subalgebra of horizontal forms for $\Omega(P)$, and such
that $D=D_\omega$. The connection $\omega$ corresponds to the
Levi-Civita connection in classical geometry.

\section{Concluding Examples $\&$ Remarks}

\noindent{\it Odd-dimensional Classes}

\medskip
As we have mentioned in Section~3, the algebra $\WGi$ may have
non-trivial odd-dimensional cohomology classes. We shall here describe an
interesting example giving $H^3(\WGi)\neq\{0\}$. As first, we are going
to prove that the higher-order calculus on $G$ can be always maximally
addopted, to allow the appearence of $3$-dimensional characteristic
classes.

As we have seen, $2$-dimensional Chern-Simons forms are labeled by
closed $\{\adj,\sigma\}$-invariant elements of $\Gamma_{\inv}^{\sstar 2}$.
Let $S_3\subseteq\Gamma^\wedge$ be the ideal generated by the elements of
the form $d^\wedge(\psi)$, where $\psi\in\Gamma_{\inv}^{\wedge 2}$
satisfy $ \sigma(\psi)=\psi$ and $\adj(\psi)=\psi\otimes 1$.
By definition, $S_3$ is a graded-differential *-ideal.

\begin{lem}
The ideal $S_3$ is $\widehat{\phi}$-invariant, in other words
\begin{equation}
\widehat{\phi}(S_3)\subseteq S_3\grten\Gamma^\wedge+
\Gamma^\wedge\grten S_3.
\end{equation}
\end{lem}

Therefore we can pass to the factor-calculus
$\Gamma^{\3}=\Gamma^\wedge/S_3$, in the framework of which all
$\{\adj,\sigma\}$-invariant elements of $\Gamma_{\inv}^{\wedge 2}$ are
closed.

As a concrete illustration, let us assume that $G$ is a
compact simple simply connected Lie group.
The elements of $G$ are naturally interpretable as characters of
$\cal{A}$. Let us denote by $\frak{g}$ the
corresponding Lie algebra. Let us consider a calculus
$\Gamma$ over $G$ defined in the following way. We put
\begin{equation}\label{def-G}
\Gamma_{\inv}=[\frak{g}^*]^{Z(G)},
\end{equation}
with $\pi\colon\cal{A}\rightarrow\Gamma_{\inv}$ given by taking
standard differentials in points of $Z(G)$. In other words
$\cal{R}$ is consisting of elements of $\ker(\e)$ the derivatives of
which vanish on $Z(G)$.

By construction, $\Gamma$ is a bicovariant *-calculus. The
adjoint action $\adj$ is induced by the classical
$\adj$ on $\frak{g}$. The *-structure is a combination of the classical
$*$ and the inverse $g\leftrightarrow g^{-1}$, acting on the coordinate
set $Z(G)$. There exists a natural left action of
$Z(G)$ on $\Gamma_{\inv}$, given by left translations on $Z(G)$.
Furthermore, the $\circ$-structure is
given by
\begin{equation}
[\psi\circ a]_g=g(a)[\psi]_g,
\end{equation}
where $g\in Z(G)$. It follows that the braiding $\sigma$ coincides with
the standard transposition, in particular $\Gamma^\vee$ is the standard
exterior algebra. However $\Gamma^\wedge\neq\Gamma^\vee$. Let us compute
the space $S_{\inv}^{\wedge 2}$. Acting by $(\pi\otimes\pi)\phi$ on the
elements of $\cal{R}$ we conclude that
\begin{lem} The universal envelope is given by the relations
\begin{equation}
S_{\inv}^{\wedge 2}=\ker(I-\sigma)\cap\Bigl\{\psi\in\Gamma^{\otimes 2}_{
\inv}\vert (g\otimes g^{-1})\psi=\psi\quad\forall g\in Z(G)\Bigr\}.
\end{equation}
\end{lem}

In particular, there exist $\adj$-singlets in the space
$\ker(I-\sigma)/S_{\inv}^{\wedge 2}$ of $\sigma$-invariant elements of
$\Gamma^{\wedge 2}_{\inv}$. Finally, let us assume that the higher-order
calculus on $G$ is described by $\Gamma^{\3}$. Having in mind
that $\adj\colon\frak{g}^*\rightarrow\frak{g}^*\otimes \cal{A}$ is
irreducible it follows that
\begin{equation}
H^3(\WGi)=\Bigl\{\mbox{$\{\adj,\sigma\}$-invariant elements of
$\Gamma_{\inv}^{\wedge 2}$}\Bigr\}.
\end{equation}
Actually in this case all $\adj$-invariant elements will be
$\sigma$-invariant.

\medskip
\noindent{\it Some Particular Cases}

\medskip
Let us consider the trivial principal bundle $P$ over a $1$-point set
$M=\{*\}$, with the trivial differential calculus. Then
$\hbim{u}=H_u^*$, in a natural manner. The corresponding contraction
maps are given by
$$ \langle\varphi,x\rangle_u^+=\varphi(x)\qquad\langle
x,\varphi\rangle_u^-
=\varphi C_u^{-1}(x).$$
Let us consider a linear map $A\colon H_u^*\rightarrow H_v^*$. Its right
transposed coincides with the standard transposed $A^\top$, while
$$ A_\perp=C_uA^\top C_v^{-1}.$$
Hence, $A$ will be transposable iff $A^\top$ intertwines $C_u$ and
$C_v$. In particular, $\hbim{u}\sim\hbim{v}$ iff $C_u$ and $C_v$ have
the same eigenvalues with the same multiplicities.

The formalism includes non-trivial examples with differential
calculus over such trivial spaces. However, it may be necessary to
refine the above introduced equivalence relation on $\RepG$,
taking the presence of $\Omega(M)$ into account.

For example, we can choose
$\Omega(P)=\Gamma^\wedge$ and assume that the calculus on $G$ is
described by $\Gamma^\vee$. In this case we have
$$ \widehat{F}=(\id\otimes\between)\widehat{\phi}.$$
Actually, the algebra
$\Omega(M)$ measures in a certain sense the difference between algebras
$\Gamma^\wedge$ and $\Gamma^\vee$. It follows that
$$\Omega(M)=\Bbb{C}\quad\Leftrightarrow\quad\Gamma^\wedge=\Gamma^\vee.$$

If $P$ is a line bundle and if the calculus on $G$ is $1$-dimensional
then the spectral sequence always stabilizes at the third term. Hence,
it is possible to pass to a long exact sequence of cohomology groups
\begin{equation*}
@>>>H^n(P)@>>>H^{n-1}(M)@>>>H^{n+1}(M)@>>>H^{n+1}(P)@>>>
\end{equation*}
as in classical geometry. In particular, if $H(P)=\Bbb{C}$ then $H(M)$
will be freely generated by the Euler class.

\medskip
\noindent{\it Non-Standard Chern Classes}

\medskip
A simple example of the appearence of infinitely many universal Chern
classes is given by $G=S_\mu U(2)$, with a $4D$-calculus $\Gamma$. Let
us consider this in more details. We have
\begin{equation*}
u=\begin{pmatrix}\alpha&-\mu\gamma^*\\
\gamma&\phantom{-\mu}\alpha^*\end{pmatrix}\qquad
C_u=\begin{pmatrix}\mu^{-1}&0\\0&
\mu^{\phantom{-1}}\end{pmatrix}
\end{equation*}
and the space $\Gamma_{\inv}$ is spanned by elements
$$
\tau=\pi(\mu^2\alpha+\alpha^*)\qquad \eta_3=
\pi(\alpha-\alpha^*)\qquad\eta_-=\pi(\gamma^*)\qquad\eta_+=\pi(\gamma),
$$
forming $\adj$-singlet and triplet respectively. Using explicit formulas
\cite{D1} for the braiding $\sigma$ it follows that $\S$ is
generated by the relations
$$\tau\eta_i=\eta_i\tau=-\frac{(1-\mu^3)\varkappa_i}{(1-\mu^2)(1+\mu)},$$
where
\begin{gather*}
\varkappa_+=\eta_+\otimes\eta_3-\mu^2\eta_3\otimes\eta_+\qquad
\varkappa_-=\eta_3\otimes\eta_--\mu^2\eta_-\otimes\eta_3\\
\varkappa_3=(1-\mu^2)\eta_3\otimes\eta_3+\mu(1+\mu^2)(\eta_+
\otimes\eta_--\eta_-\otimes\eta_+).
\end{gather*}

Let $\S_3\subseteq\S$ be the subalgebra generated by triplet
elements $\eta_i$. The above formulas imply that $\S_3$ is
generated by the following cubic relations
$$ \eta_i\varkappa_j=\varkappa_i\eta_j,\qquad i,j\in\{{+},3,{-}\}.$$

Every element $w\in\S$ can be uniquely decomposed in the form
$$w=p(\tau)+q(\eta_+,\eta_3,\eta_-),$$
where $p,q$ are polynomial expressions. The simplest way to prove the
non-triviality of Chern classes in high orders is to factorize $\S$
through the ideal generated by $\{\eta_+,\eta_-\}$. The factor-algebra
has two generators $\tau,\eta_3$, satisfying a single relation
$$\eta_3\tau=\tau\eta_3=-\frac{1-\mu^3}{1+\mu^{\phantom{3}}}\eta_3^2.$$
A direct calculation gives
\begin{multline*}
\Delta^\lambda(\chi_u)=\\
=\Bigl(1+\frac{\tau\lambda}{1+\mu^2}\Bigr)^{\mu
+1/\mu}+\Bigl(1+\frac{\mu\lambda\eta_3}{1+\mu}\Bigr)^\mu
\Bigl(1-\frac{\eta_3\lambda}{1+\mu}\Bigr)^{1/\mu}
-\Bigl(1-\frac{(1-\mu^3)\lambda\eta_3}{(1+\mu)(1+\mu^2)}\Bigr)^{\mu+
1/\mu}
\end{multline*}
for the factorized $\Delta$.

\medskip
\noindent{\it Equivalence of Vector Bundles}

\medskip

As we have already mentioned, the introduced equivalence relation
$\sim$ between vector bundles $\hbim{u}$ depends also on the
specification of the calculus (via the algebra $\hor(P)$). The
additional structural element coming from the calculus is expressed via
canonical flip-over operators $\sigma_u\colon\bim{u}\otimes_M\!\Omega(M)
\leftrightarrow\Omega(M)\otimes_M\!\bim{u}$. Homomorphisms acting on pure
sections $\bim{u}$ will be extendible to the higher-order level iff
they commute with these flip-over operators.

Covariant derivative maps $D\in\rhom(\hbim{u})$ are in a natural
correspondence with their restrictions $D\colon\bim{u}\rightarrow
\bim{u}\otimes_M\!\Omega^1(M)$, characterized as first-order linear maps
satisfying
\begin{equation*}
D(\psi f)=D(\psi)f+\psi\otimes f,
\end{equation*}
for each $\psi\in\bim{u}$ and $f\in\cal{V}$. In other words, the
maps $D$ are {\it connections} on a finite projective
module $\bim{u}$, in the
sense of \cite{C1}. A given connection $D$ will satisfy the left Leibniz
rule iff the following compatibility condition holds
$$[\sigma_uD\sigma_u^{-1}](\alpha\otimes\psi)=d\alpha\otimes\psi+
(-1)^{\partial\alpha}\alpha\sigma_uD(\psi),$$
where $\alpha\in\Omega(M)$ and $\psi\in\bim{u}$.


\newpage
\begin{thebibliography}{10}
\bibitem{C1} Connes A: {\it Non-commutative differential geometry},
Extrait des Publications Math\'ematiques--IHES {\bf 62} (1986)
\bibitem{C2} Connes A: {\it Geometrie non commutative}, InterEditions
Paris (1990)
\bibitem{bot} Bott R $\&$ Tu W: {\it Differential Forms in Algebraic
Topology}, Springer-Verlag New-York (1982)
\bibitem{b}  Brzezi\'nski T $\&$ Majid S: {\it Quantum Group Gauge
Theory on Quantum Spaces}, CMP {\bf 157} 591--638 (1993)
\bibitem{D1} \mbox{{\Dj}ur{\dj}evi\'c M}
{\it Geometry of Quantum Principal Bundles I}, Preprint QmmP 6/92,
Belgrade University, Serbia (To appear in CMP)
\bibitem{D2} \mbox{{\Dj}ur{\dj}evi\'c M}:
{\it Geometry of Quantum Principal Bundles II--Extended Version},
Preprint, Instituto de Matematicas, UNAM, M\'exico (1994)
\bibitem{D3} \mbox{{\Dj}ur{\dj}evi\'c M}:
{\it Quantum Principal Bundles $\&$ Tannaka-Krein Duality Theory},
Preprint, Instituto de Matematicas, UNAM, M\'exico (1995)
\bibitem{diff} \mbox{{\Dj}ur{\dj}evi\'c M}:
{\it On Differential Structures on Quantum Principal Bundles},
Preprint, Instituto de Matematicas, UNAM, M\'exico (1994)
\bibitem{frm1} \mbox{{\Dj}ur{\dj}evi\'c M}:
{\it On Framed Quantum Principal Bundles},
Preprint, Instituto de Matematicas, UNAM, M\'exico (1994)
\bibitem{frm2} \mbox{{\Dj}ur{\dj}evi\'c M}:
{\it General Frame Structures on Quantum Principal Bundles},
Preprint, Instituto de Matematicas, UNAM, M\'exico (1995)
\bibitem{KN} Kobayashi S and Nomizu K: {\it Foundations of Differential
Geometry}, Interscience Publishers New York, London (1963)
\bibitem{pfl} Pflaum M: {\it Quantum Groups on Fibre Bundles},
CMP (1994)
\bibitem{W1} Woronowicz S L: {\it Compact Matrix Pseudogroups}, CMP
{\bf 111} 613--665 (1987)
\bibitem{W2} Woronowicz S L: {\it Differential Calculus on Compact
Matrix Pseudogroups/ Quantum Groups}, CMP {\bf 122} 125--170 (1989)
\end{thebibliography}
\end{document}